\documentclass[10pt,english]{amsart}
\usepackage{helvet}
\usepackage[T1]{fontenc}
\usepackage[latin1]{inputenc}
\usepackage{a4wide}
\usepackage{amssymb}
\usepackage{babel}

 \theoremstyle{plain}
 \newtheorem{thm}{Theorem}[section]
 \numberwithin{equation}{section}  
 \numberwithin{figure}{section} 
 \theoremstyle{plain}
 \theoremstyle{plain}
 \newtheorem*{thm*}{Theorem}
 \theoremstyle{definition}
 \newtheorem{defn}[thm]{Definition}
 \theoremstyle{remark}
 \newtheorem{rem}[thm]{Remark}
 \theoremstyle{definition}
  \newtheorem{example}[thm]{Example}
\theoremstyle{definition}
  \newtheorem{notation}[thm]{Notation}
 \theoremstyle{plain}
 \newtheorem{lem}[thm]{Lemma}  
 \theoremstyle{plain}
 \newtheorem{prop}[thm]{Proposition}  
 \theoremstyle{plain}
 \newtheorem{cor}[thm]{Corollary} 
 \theoremstyle{definition}
  \newtheorem*{example*}{Example}

\usepackage{latexsym,amsfonts, xypic,amssymb}
\usepackage{pstricks,pst-node,pst-tree}
\xyoption{all}

\theoremstyle{remark}

\theoremstyle{definition}
\newtheorem{enavant}[thm]{}

\begin{document}

\title{On a class of Danielewski surfaces in affine $3$-space }

\author{Adrien Dubouloz}
\address{Institut Fourier, UMR 5582 CNRS-UJF, BP 74, 38402 Saint-Martin d'H\`eres
Cedex, France}
\email{adrien.dubouloz@ujf-grenoble.fr}

\author{Pierre-Marie Poloni}
\address{Institut de Math\'ematiques de Bourgogne, CNRS-UMR 5584, 9 avenue Alain
Savary, BP 47870, 21078 Dijon, France}
\email{pierre-marie.poloni@u-bourgogne.fr }

\maketitle
 \insert\footins{\footnotesize\noindent \textbf{Mathematics Subject Classification (2000)}: 14R10, 14R05.


 \unskip\strut\par}

\begin{abstract}
In \cite{ML90} and \cite{ML01}, L. Makar-Limanov computed the automorphism
groups of surfaces in $\mathbb{C}^{3}$ defined by the equations $x^{n}z-P\left(y\right)=0$,
where $n\geq1$ and $P\left(y\right)$ is a nonzero polynomial. Similar
results have been obtained by A. Crachiola \cite{Cra06} for surfaces
defined by the equations $x^{n}z-y^{2}-\sigma\left(x\right)y=0$, where $n\geq2$
and $\sigma\left(0\right)\neq0$, defined over an arbitrary base field.
Here we consider the more general surfaces defined by the equations $x^{n}z-Q\left(x,y\right)=0$,
where $n\geq2$ and $Q\left(x,y\right)$ is a polynomial with coefficients
in an arbitrary base field $k$. We characterise among them
the ones which are Danielewski surfaces in the sense of \cite{DubG03},
and we compute their automorphism groups. We study closed embeddings
of these surfaces in affine $3$-space. We show that in general their
automorphisms do not extend to the ambient space. Finally, we give
explicit examples of $\mathbb{C}^{*}$-actions on a surface in $\mathbb{A}_{\mathbb{C}}^{3}$
which can be extended holomorphically but not algebraically to  $\mathbb{C}^{*}$-actions
on $\mathbb{A}_{\mathbb{C}}^{3}$.
\end{abstract}

\section*{Introduction }

Since they appeared in a celebrated counterexample to the
Cancellation Problem due to W. Danielewski \cite{Dan89}, the surfaces
defined by the equations $xz-y\left(y-1\right)=0$ and $x^{2}z-y\left(y-1\right)=0$
in $\mathbb{C}^{3}$ and their natural generalisations, such as surfaces
defined by the equations $x^{n}z-P\left(y\right)=0$, where $P\left(y\right)$
is a nonconstant polynomial, have been studied in many different
contexts. Of particular interest is the fact
that they can be equipped with nontrivial actions of the additive
group $\mathbb{C}_{+}$. The general orbits of these actions coincide with the general fibers of  $\mathbb{A}^{1}$-fibrations $\pi:S\rightarrow\mathbb{A}^{1}$,
that is,  surjective morphisms with generic fiber isomorphic to an
affine line. Normal affine
surfaces $S$ equipped with an $\mathbb{A}^{1}$-fibration $\pi:S\rightarrow\mathbb{A}^{1}$
can be roughly classified into two classes according the following
alternative : either $\pi:S\rightarrow\mathbb{A}^{1}$ is a unique
$\mathbb{A}^{1}$-fibration on $S$ up to automorphisms of the base,
or there exists a second $\mathbb{A}^{1}$-fibration $\pi':S\rightarrow\mathbb{A}^{1}$
with general fibers distinct from the ones of $\pi$.

Due to the symmetry between the variables $x$ and $z$, a surface
defined by the equation $xz-P\left(y\right)=0$ admits two distinct $\mathbb{A}^1$-fibrations
over the affine line. In contrast, it was established by L. Makar-Limanov \cite{ML01} that
on a surface $S_{P,n}$ defined by the equation $x^{n}z-P\left(y\right)=0$
in $\mathbb{C}^{3}$, where $n\geq2$ and where $P\left(y\right)$
is a polynomial of degree $r\geq2$, the projection ${\rm pr}_{x}:S_{P,n}\rightarrow\mathbb{C}$
is a unique $\mathbb{A}^{1}$-fibration up to automorphisms of the
base. In his proof, L. Makar-Limanov used the correspondence between algebraic $\mathbb{C}_+$-actions on an affine surface $S$ and locally nilpotent derivations of the algebra of regular functions on $S$. It turns out that his proof is essentially independent of the base field $k$ provided that we replace locally nilpotent derivations by suitable systems of Hasse-Schmidt derivations when the characteristic of $k$ is positive (see e.g., \cite{Cra06}).

The fact that  an affine surface $S$ admits a unique $\mathbb{A}^{1}$-fibration
$\pi:S\rightarrow\mathbb{A}^{1}$ makes its study simpler. For instance,
every automorphism of $S$ must preserve this fibration. In this context,
a result due to J. Bertin \cite{Ber83} asserts that the identity
component of the automorphisms group of such a surface is an algebraic
pro-group obtained as an increasing union of solvable algebraic subgroups
of rank $\leq1$. For surfaces defined by the equations $x^{n}z-P\left(y\right)=0$
in $\mathbb{C}^{3}$, the picture has been completed by L. Makar-Limanov
\cite{ML01} who gave explicit generators of their automorphisms
groups. Similar results have been obtained over arbitrary base fields by A. Crachiola
\cite{Cra06} for surfaces defined by the equations $x^{n}z-y^{2}-\sigma\left(x\right)y=0$, where $\sigma\left(x\right)$ is a polynomial such that $\sigma\left(0\right)\neq 0$.

The latter surfaces are particular examples of a general class of
$\mathbb{A}^{1}$-fibred surfaces called \emph{Danielewski surfaces} \cite{DubG03}, that is, normal integral affine surface $S$ equipped with an $\mathbb{A}^{1}$-fibration
$\pi:S\rightarrow\mathbb{A}_{k}^{1}$ over an affine line with a fixed $k$-rational point
$o$, such that every fiber $\pi^{-1}\left(x\right)$, where $x\in\mathbb{A}_{k}^{1}\setminus\left\{ o\right\} $, is geometrically integral, and such that every irreducible component
of $\pi^{-1}\left(o\right)$ is geometrically integral. In this article, we consider Danielewski surfaces  $S_{Q,n}$ in $\mathbb{A}_{k}^{3}$ defined by an equation of the form $x^{n}z-Q\left(x,y\right)=0$, where $n\geq2$ and where $Q\left(x,y\right)\in k\left[x,y\right]$
is a polynomial such that $Q\left(0,y\right)$ splits with $r\geq2$
simple roots in $k$. This class contains most of the surfaces considered by L. Makar-Limanov and A. Crachiola.

The paper is organised as follows. First, we briefly recall definitions about weighted rooted trees and  the notion of equivalence of algebraic surfaces in an affine $3$-space.
In section $2$, we recall from \cite{DubG03} the
main facts about Danielewski surfaces and we review the correspondence between these surfaces and certain classes of  weighted trees in a form appropriate to our needs. We also
generalise to arbitrary base fields $k$ some results which are only
stated for fields of characteristic zero in \cite{Dub02} and \cite{DubG03}.
In particular, the case of Danielewski surfaces which admit two $\mathbb{A}^{1}$-fibrations with distinct general
fibers is studied in  Theorem \ref{thm:Comb_ML_Trivial}. We show that these surfaces correspond to Danielewski surfaces $S\left(\gamma\right)$
defined by the fine $k$-weighted trees $\gamma$ which are called \emph{combs} and we give explicit embeddings of them. This result generalises Theorem 4.2 in \cite{DubEmb05}. 

In section $3$, we classify Danielewski surfaces $S_{Q,h}$ in $\mathbb{A}_{k}^{3}$
defined by equations of the form $x^{h}z-Q\left(x,y\right)=0$ and determine their automorphism groups. We remark that such a surface admits many embeddings as a surface $S_{Q,h}$. In particular, we establish in Theorem \ref{thm:Equivalent_charac} that these surfaces can always be embedded as surface $S_{\sigma,h}$ defined by an equation of the form  $x^{h}z-\prod_{i=1}^r\left(y-\sigma_i\left(x\right)\right)=0$ for a suitable collection of polynomials $\sigma=\left\{\sigma_i\left(x\right)\right\}_{i=1,\ldots,r}$. We say that these surfaces $S_{\sigma,h}$ are \emph{standard form} of Danielewski surfaces $S_{Q,h}$. 
Next, we compute ( Theorem \ref{thm:Main_auto_thm}) the automorphism groups of Danielewski surfaces in standard form. We show in particular that any of them comes as the restriction of an algebraic automorphism of the ambient space.

Finally we consider the problem of extending  automorphisms of a given Danielewski surface $S_{Q,h}$ to automorphisms of the ambient space $\mathbb{A}^3_k$. We show that this is always possible in the holomorphic category but not in the algebraic one. We give explicit examples which come from the study of multiplicative group actions on Danielewski surfaces. For instance, we prove that
every surface $S\subset\mathbb{A}_{\mathbb{C}}^{3}$ defined by the equation $x^{h}z-\left(1-x\right)P\left(y\right)=0$, where $h\geq2$ and where
$P\left(y\right)$ has $r\geq2$ simple roots, admits a nontrivial
$\mathbb{C}^{*}$-action which is algebraically inextendable but holomorphically
extendable to $\mathbb{A}_{\mathbb{C}}^{3}$. In particular, even the involution of the surface $S$ defined by the equation $x^{2}z-\left(1-x\right)P\left(y\right)=0$ induced by the endomorphism
$J\left(x,y,z\right)=\left(-x,y,\left(1+x\right)\left(\left(1+x\right)z+P\left(y\right)\right)\right)$
of $\mathbb{A}^3_k$ does not extend to an algebraic automorphism of $\mathbb{A}_{k}^{3}$.

\section{Preliminaries}

\subsection{Basic facts on weighted rooted trees}
\begin{defn}
A \emph{tree} is a nonempty, finite, partially ordered set $\Gamma=\left(\Gamma,\leq\right)$
with a unique minimal element $e_{0}$ called the \emph{root}, and
such that for every $e\in\Gamma$ the subset $\left(\downarrow e\right)_{\Gamma}=\left\{ e'\in\Gamma,\, e'\leq e\right\} $
is a chain \emph{}for the induced ordering.
\end{defn}
\begin{enavant} \label{txt-def:subchain_def} A minimal sub-chain
$\overleftarrow{e'e}=\left\{ e'<e\right\} $ with two elements of
a tree $\Gamma$ is called \emph{an edge} of $\Gamma$. We denote
the set of all edges in $\Gamma$ by $E\left(\Gamma\right)$. An element
$e\in\Gamma$ such that $\textrm{Card}\left(\downarrow e\right)_{\Gamma}=m$
is said to be \emph{at level} $m$. The maximal elements $e_{i}=e_{i,m_{i}}$,
where $m_{i}=\textrm{Card}\left(\downarrow e_{i}\right)_{\Gamma}$
of $\Gamma$ are called the \emph{leaves} of $\Gamma$. We denote
the set of those elements by $L\left(\Gamma\right)$. The maximal
chains of $\Gamma$ are the chains \begin{equation}
\Gamma_{e_{i,m_{i}}}=\left(\downarrow e_{i,m_{i}}\right)_{\Gamma}=\left\{ e_{i,0}=e_{0}<e_{i,1}<\cdots<e_{i,m_{i}}\right\} ,\quad e_{i,m_{i}}\in L\left(\Gamma\right).\label{eq:Maximal_Chain_Notation}\end{equation}
 We say that $\Gamma$ has \emph{height} $h=\max\left(m_{i}\right)$.
The \emph{children} of an element $e\in\Gamma$ are the elements of
$\Gamma$ at relative level $1$ with respect to $e$, i.e.,
the maximal elements of the subset $\left\{ e'\in\Gamma,\: e'>e\right\} $
of $\Gamma$.

\end{enavant}

\begin{defn}
\label{def:fine_weighted_def1} A \emph{fine} $k$\emph{-weighted
tree} $\gamma=\left(\Gamma,w\right)$ is a tree $\Gamma$ equipped
with a weight function $w:E\left(\Gamma\right)\rightarrow k$ with
values in a field $k$, which assigns an element $w\left(\overleftarrow{e'e}\right)$
of $k$ to every edge $\overleftarrow{e'e}$ of $\Gamma$, in such
a way that $w\left(\overleftarrow{e'e_{1}}\right)\neq w\left(\overleftarrow{e'e_{2}}\right)$
whenever $e_{1}$ and $e_{2}$ are distinct children of a same element $e'$.
\end{defn}
\noindent In what follows, we frequently consider the following classes
of trees.

\begin{defn} \label{CombDef} Let $\Gamma$ be a rooted tree.

a) If all the leaves of $\Gamma$ are at the same level $h\geq 1$ and if there exists a unique element $\bar{e}_{0}\in\Gamma$ for which $\Gamma\setminus\left\{ \bar{e}_{0}\right\} $ is a nonempty disjoint union of chains then we say that $\Gamma$ is a \emph{rake}.

b) If $\Gamma\setminus L\left(\Gamma\right)$ is a chain then we say that $\Gamma$ is a \emph{comb}. Equivalently, $\Gamma$ is a comb if and only if every $e\in\Gamma\setminus L\left(\Gamma\right)$
has at most one child which is not a leaf of $\Gamma$.

\end{defn}
\begin{figure}[h]

\begin{pspicture}(1,0.6)(8,-1.5)

\rput(1,0){

\pstree[treemode=R,radius=2.5pt,treesep=0.5cm,levelsep=0.6cm]{\Tc{3pt}~[tnpos=a]{$e_0$}}{

\pstree{\TC*}{\skiplevels{1}\pstree{\TC*~[tnpos=a]{$\bar{e}_0$}}{

  \pstree{\TC*}{\skiplevels{2} \TC*\endskiplevels}

 \pstree{\TC*}{\skiplevels{2} \TC*\endskiplevels}

 \pstree{\TC*}{\skiplevels{2} \TC*\endskiplevels}

}\endskiplevels}

}

\rput(-1.5,-1.2){ A rake rooted in $e_0$.}

\rput(5,0){

\pstree[treemode=R,radius=2.5pt,treesep=0.5cm,levelsep=1cm]{\Tc{3pt}~[tnpos=a]{$e_0$}}{

 \pstree{\TC*}{  \pstree{\TC*}   {\Tn \pstree{\TC*} {\TC*\TC*\TC*} \TC* }   }

}

}

}

\rput(8,-1.2){ A comb rooted in $e_0$.}
\end{pspicture}

\end{figure}

\subsection{Algebraic and analytic equivalence of closed embeddings}

\indent\newline\noindent Here we briefly discuss the notions of algebraic
and analytic equivalences of closed embeddings of a given affine algebraic
surface in an affine $3$-space.

Let $S$ be an irreducible affine surface and let $i_{P_1}:S\hookrightarrow \mathbb{A}^3_{k}$ and $i_{P_2}:S\hookrightarrow \mathbb{A}^3_{k}$ be embeddings of $S$ in a same affine $3$-space  as closed subschemes defined by  polynomial equations $P_{1}=0$ and $P_{2}=0$ respectively.

\begin{defn}
\label{def:Algebraic_equiv_def} In the above setting, we say that the closed embeddings $i_{P_1}$ and $i_{P_2}$  are \emph{algebraically equivalent}
if one of the following equivalent conditions is satisfied:

1) There exists an automorphism $\Phi$ of $\mathbb{A}^3_{k}$ such that $i_{P_2}=i_{P_1}\circ\Phi$.

2) There exists an automorphism $\Phi$ of $\mathbb{A}_{k}^{3}$ and a nonzero constant
$\lambda\in k^{*}$ such that  $\Phi^{*}P_{1}=\lambda P_{2}$.

3) There exists automorphisms $\Phi$ and $\phi$  of $\mathbb{A}_{k}^{3}$
and $\mathbb{A}_{k}^{1}$ respectively such that $P_{2}\circ\Phi=\phi\circ P_{1}$.
\end{defn}

\begin{enavant} \label{txt:Analytic_equivalence_def} Over the field $k=\mathbb{C}$ of complex numbers, one can also consider holomorphic automorphisms. With the notation of definition \ref{def:Algebraic_equiv_def}, two closed algebraic embeddings $i_{P_1}$ and $i_{P_2}$  of a given affine surface $S$ in $\mathbb{A}^3_{\mathbb{C}}$ are called \emph{holomorphically equivalent} if there exists a biholomorphism $\Phi:\mathbb{A}_{\mathbb{C}}^{3}\rightarrow\mathbb{A}_{\mathbb{C}}^{3}$ such that $i_{P_2}=i_{P_1}\circ\Phi$. Clearly, the embeddings  $i_{P_2}$ and $i_{P_1}$ are holomorphically equivalent if and only if there exists a biholomorphism $\Phi:\mathbb{A}_{\mathbb{C}}^{3}\rightarrow\mathbb{A}_{\mathbb{C}}^{3}$ such that $\Phi^{*}\left(P_1\right)=\lambda P_2$ for a certain nowhere vanishing holomorphic function $\lambda$. Since there are many nonconstant holomorphic functions with this property on $\mathbb{A}_{\mathbb{C}}^{3}$, $\Phi$ need not preserve the algebraic families of level surfaces $P_{1}:\mathbb{A}_{\mathbb{C}}^{3}\rightarrow\mathbb{A}_{\mathbb{C}}^{1}$
and $P_{2}:\mathbb{A}_{\mathbb{C}}^{n}\rightarrow\mathbb{A}_{\mathbb{C}}^{1}$. So holomorphic equivalence  is a weaker requirement than algebraic equivalence.
\end{enavant}

\section{Danielewski surfaces }

For certain authors, a Danielewski surface
is an affine surface $S$ which is algebraically isomorphic to a surface
in $\mathbb{C}^{3}$ defined by an equation of the form $x^{n}z-P\left(y\right)=0$,
where $n\geq1$ and $P\left(y\right)\in\mathbb{C}\left[y\right]$. These surfaces come
equipped with a surjective morphism $\pi={\rm pr}_{x}\mid_{S}:S\rightarrow\mathbb{A}^{1}$
restricting to a trivial $\mathbb{A}^{1}$-bundle over the complement
of the origin. Moreover, if the roots $y_{1},\ldots,y_{r}\in\mathbb{C}$
of $P\left(y\right)$ are simple, then the fibration $\pi={\rm pr}_{x}\mid_{S}:S\rightarrow\mathbb{A}^{1}$
factors through a locally trivial fiber bundle over the affine line
with an $r$ -fold origin (see e.g., \cite{Dan89} and \cite{Fie94}).
In \cite{DubG03}, the first author used the term Danielewski surface
to refer to an affine surface $S$ equipped with a morphism $\pi:S\rightarrow\mathbb{A}^{1}$
which factors through a locally trivial fiber bundle in a similar
way as above. In what follows, we keep this point of view, which leads
to a natural geometric generalisation of the surfaces constructed
by W. Danielewski \cite{Dan89}. We recall that an $\mathbb{A}^{1}$\emph{-fibration}
over an integral scheme $Y$ is a faithfully flat (i.e., flat and
surjective) affine morphism $\pi:X\rightarrow Y$ with generic fiber
isomorphic to the affine line $\mathbb{A}_{K\left(Y\right)}^{1}$
over the function field $K\left(Y\right)$ of $Y$. The following
definition is a generalisation to arbitrary base fields $k$ of the
one introduced in \cite{DubG03}.

\begin{defn}
\label{def:DanSurf_Def} A \emph{Danielewski surface} is an integral
affine surface $S$ defined over a field $k$, equipped with an $\mathbb{A}^{1}$-fibration
$\pi:S\rightarrow\mathbb{A}_{k}^{1}$ restricting to a trivial $\mathbb{A}^{1}$-bundle
over the complement of the a $k$-rational point $o$ of $\mathbb{A}_{k}^{1}$ and
such that the fiber $\pi^{-1}\left(o\right)$ is reduced, consisting
of a disjoint union of affine lines $\mathbb{A}_{k}^{1}$ over $k$.
\end{defn}

\begin{notation} In what follows, we fix an isomorphism $\mathbb{A}^1_k\simeq {\rm Spec}\left(k\left[x\right]\right)$ and we assume that the $k$-rational point  $o$  is simply the "origin" of $\mathbb{A}^1_k$, that is, the closed point $(x)$ of ${\rm Spec}\left(k\left[x\right]\right)$.
\end{notation}

\begin{enavant} In the following subsections, we recall the correspondence between Danielewski surfaces and weighted rooted trees established by the first author in \cite{DubG03} in a form
appropriate to our needs. Although the results given in \emph{loc}.
\emph{cit}. are formulated for surfaces defined over a field of characteristic
zero, most of them remain valid without any changes over a field of
arbitrary characteristic. We provide full proofs only when additional
arguments are needed. Then we consider Danielewski surfaces
$S$ with a trivial canonical sheaf $\omega_{S/k}=\Lambda^{2}\Omega_{S/k}^{1}$.
We call them  \emph{special Danielewski surfaces}. We give a complete
classification of these surfaces in terms of their associated weighted
trees.

\end{enavant}

\subsection{Danielewski surfaces and weighted trees}

\indent\newline\noindent Here we review the correspondence  which associates  to every fine $k$-weighted tree $\gamma=\left(\Gamma,w\right)$
a Danielewski surface $\pi:S\left(\gamma\right)\rightarrow\mathbb{A}_{k}^{1}={\rm Spec}\left(k\left[x\right]\right)$
which is the total space of an $\mathbb{A}^{1}$-bundle over
the scheme $\delta:X\left(r\right)\rightarrow\mathbb{A}_{k}^{1}$
obtained from $\mathbb{A}_{k}^{1}$ by replacing its origin $o$ by
$r\geq1$ $k$-rational points  $o_{1},\ldots,o_{r}$.

\begin{notation}
In what follows we denote
by $\mathcal{U}_{r}=\left(X_{i}\left(r\right)\right)_{i=1,\ldots,r}$
the canonical open covering of $X\left(r\right)$ by means of the
subsets $X_{i}\left(r\right)=\delta^{-1}\left(\mathbb{A}_{k}^{1}\setminus\left\{ o\right\} \right)\cup\left\{ o_{i}\right\} \simeq\mathbb{A}_{k}^{1}$.
\end{notation}

\begin{enavant} \label{txt:Abstract_DS_morph}\label{pro:WeightedTree_2_DanSurf}
Let $\gamma=\left(\Gamma,w\right)$ be a fine $k$-weighted tree $\gamma=\left(\Gamma,w\right)$ of
height $h$, with leaves $e_{i}$ at levels $n_{i}\leq$$h$, $i=1,\ldots,r$. To every
 maximal sub-chain $\gamma_{i}=\left(\downarrow e_{i}\right)$
of $\gamma$ (see \ref{txt-def:subchain_def} for the notation) we associate a polynomial
\[
\sigma_{i}\left(x\right)=\sum_{j=0}^{n_{i}-1}w\left(\overleftarrow{e_{i,j}e_{i,j+1}}\right)x^{j}\in k\left[x\right],\quad i=1,\ldots,r.\]
 We let $\rho:S\left(\gamma\right)\rightarrow X\left(r\right)$ be
the unique $\mathbb{A}^{1}$-bundle over $X\left(r\right)$ which
becomes trivial on the canonical open covering $\mathcal{U}_{r}$, and is
defined by pairs of transition functions \[
\left(f_{ij},g_{ij}\right)=\left(x^{n_{j}-n_{i}},x^{-n_{i}}\left(\sigma_{j}\left(x\right)-\sigma_{i}\left(x\right)\right)\right)\in k\left[x,x^{-1}\right]^{2},\quad i,j=1,\ldots,r.\]
 This means that $S\left(\gamma\right)$ is obtained by gluing $n$
copies $S_{i}=\textrm{Spec}\left(k\left[x\right]\left[u_{i}\right]\right)$
of the affine plane $\mathbb{A}_{k}^{2}$ over $\mathbb{A}_{k}^{1}\setminus\left\{o\right\} \simeq\textrm{Spec}\left(k\left[x,x^{-1}\right]\right)$
by means of the transition isomorphisms induced by the $k\left[x,x^{-1}\right]$-algebras
isomorphisms\[
k\left[x,x^{-1}\right]\left[u_{i}\right]\stackrel{\sim}{\rightarrow}k\left[x,x^{-1}\right]\left[u_{j}\right],\quad u_{i}\mapsto x^{n_{j}-n_{i}}u_{j}+x^{-n_{i}}\left(\sigma_{j}\left(x\right)-\sigma_{i}\left(x\right)\right)\qquad i\neq,\, i,j=1,\ldots,r.\]
 This definition makes sense as the transition functions
$g_{ij}$ satisfy the twisted cocycle relation $g_{ik}=g_{ij}+x^{n_{j}-n_{i}}g_{jk}$
in $k\left[x,x^{-1}\right]$ for every triple of distinct indices $i$, $j$ and $k$. Since $\gamma$ is a fine weighted tree, it follows that for every pair of distinct indices $i$ and $j$, the rational function  $g_{ij}=x^{-n_i}\left(\sigma_j\left(x\right)-\sigma_i\left(x\right)\right)\in k\left[x,x^{-1}\right]$ does not extend to a regular function on $\mathbb{A}^1_k$. This condition guarantees that
$S\left(\gamma\right)$ is a separated scheme, whence an affine surface by
virtue of Fieseler's criterion (see proposition 1.4 in \cite{Fie94}).
Therefore, $\pi_{\gamma}=\delta\circ\rho:S\left(\gamma\right)\rightarrow\mathbb{A}_{k}^{1}=\textrm{Spec}\left(k\left[x\right]\right)$
is a Danielewski surface, the fiber $\pi^{-1}\left(o\right)$ being
the disjoint union of affine lines \[
C_{i}=\pi_{\gamma}^{-1}\left(o\right)\cap S_{i}\simeq\textrm{Spec}\left(k\left[u_{i}\right]\right),\quad i=1,\ldots,r.\]
 \end{enavant}

\begin{enavant} \label{txt:canonical_morphism} A Danielewski surface $\pi:S\left(\gamma\right)\rightarrow\mathbb{A}_{k}^{1}$
above comes canonically equipped with a birational morphism $\left(\pi,\psi_{\gamma}\right):S\rightarrow\mathbb{A}_{k}^{1}\times\mathbb{A}_{k}^{1}=\textrm{Spec}\left(k\left[x\right]\left[t\right]\right)$
restricting to an isomorphism over $\mathbb{A}_{k}^{1}\setminus\left\{ o\right\} $. Indeed, this morphism corresponds  to the unique regular function $\psi_{\gamma}$ on $S\left(\gamma\right)$
whose restrictions to the open subsets $S_i \simeq\textrm{Spec}\left(k\left[x\right]\left[u_{i}\right]\right)$ of $S$ are given by the polynomials \[\psi_{\gamma,i}=x^{n_i}u_{i}+\sigma_{i}\left(x\right)\in k\left[x\right]\left[u_i\right],\quad i=1,\ldots,r.\]
This function is referred to as the \emph{canonical
function} on $S\left(\gamma\right)$. The morphism $\left(\pi_{\gamma},\psi_{\gamma}\right):S\left(\gamma\right)\rightarrow\mathbb{A}_{k}^{2}$
is called the \emph{canonical birational morphism} from $S\left(\gamma\right)$
to $\mathbb{A}_{k}^{2}$.

\end{enavant}

\begin{enavant} It turns out that there exists a one-to-one correspondence between pairs $\left(S,\left(\pi,\psi\right)\right)$ consisting of a Danielewski surface $\pi:S\rightarrow\mathbb{A}_{k}^{1}$
and a birational morphism $\left(\pi,\psi\right):S\rightarrow\mathbb{A}_{k}^{2}$
restricting to an isomorphism outside the fiber $\pi^{-1}\left(o\right)$
and fine $k$-weighted trees $\gamma$. In particular,  Proposition
3.4 in \cite{DubG03}, which remains valid over arbitrary base fields
$k$, implies the following result.

\end{enavant}

\begin{thm}
\label{thm:GenDanMor_2_Tree} For every pair consisting of a Danielewski
surface $\pi:S\rightarrow\mathbb{A}_{k}^{1}$ and a birational morphism
$\left(\pi,\psi\right):S\rightarrow\mathbb{A}_{k}^{1}\times\mathbb{A}_{k}^{1}$
restricting to an isomorphism over $\mathbb{A}_{k}^{1}\setminus\left\{ o\right\} $,
there exists a unique fine $k$-weighted tree $\gamma$ and an isomorphism
$\phi:S\stackrel{\sim}{\rightarrow}S\left(\gamma\right)$ such that
$\psi=\psi_{\gamma}\circ\phi$.
\end{thm}
\begin{rem}
\label{rem:Psi_2_Level1} If $\gamma=\left(\Gamma,w\right)$ is not
the trivial tree with one element then the canonical function $\psi_{\gamma}:S\left(\gamma\right)\rightarrow\mathbb{A}_{k}^{1}$
on the corresponding Danielewski surface $\pi:S\left(\gamma\right)\rightarrow\mathbb{A}_{k}^{1}$
is locally constant on the fiber $\pi^{-1}\left(o\right)$. It takes
the same value on two distinct irreducible components of $\pi^{-1}\left(o\right)$
if and only if the corresponding leaves of $\gamma$ belong to a same
subtree of $\gamma$ rooted in an element at level $1$. Since every Danielewski surface nonisomorphic
to $\mathbb{A}_{k}^{2}$ admits a birational morphism $\left(\pi,\psi\right)$
for which $\psi$ is locally constant but not constant on the fiber
$\pi^{-1}\left(o\right)$, it follows that every such surface correspond to a tree $\gamma$ with at least two elements at level $1$.
\end{rem}

\subsection{$\mathbb{A}^{1}$-fibrations on Danielewski surfaces}

\indent\newline\noindent Suppose that the structural $\mathbb{A}^{1}$-fibration $\pi:S\rightarrow\mathbb{A}_{k}^{1}$  on a Danielewski surface $S$ is unique up to automorphisms of the base. Then a second Danielewski surface $\pi':S'\rightarrow\mathbb{A}_{k}^{1}$ will be isomorphic to $S$ as an abstract surface if and only if it is isomorphic to $S$  as a fibered surface, that is, if and only if there exists a commutative diagram
\[\xymatrix{ S \ar[r]^{\sim}_{\Phi} \ar[d]_{\pi} & S' \ar[d]^{\pi'} \\ \mathbb{A}^1_k \ar[r]^{\sim}_{\phi} & \mathbb{A}^1_k \;, }\]
where $\Phi:S\stackrel{\sim}{\rightarrow} S'$ is an isomorphism and $\phi$ is an automorphism of $\mathbb{A}^1_k$ preserving the origin $o$.

\begin{enavant} So it is useful to a have characterisation of those Danielewski surfaces admitting two $\mathbb{A}^1$-fibrations with distinct general fibers. The first result toward such a classification has been obtained by T. Bandman and L. Makar-Limanov \cite{BML01} who established that a complex Danielewski
surface $S$ with a trivial canonical sheaf $\omega_{S}$ admits two
$\mathbb{A}^{1}$-fibrations with distinct general fibers if and only
if it is isomorphic to a surface $S_{P,1}$ in $\mathbb{A}_{\mathbb{C}}^{3}$
defined by the equation $xz-P\left(y\right)=0$, where $P$ is a  polynomial
with simple roots. Over a field of characteristic zero, a complete classification has been given by the first author in \cite{DubG03} and \cite{DubEmb05}. It turns out that the main result of \cite{DubEmb05} remains valid over arbitrary base fields. This leads to the following characterisation.

\end{enavant}

\begin{thm}
\label{thm:Comb_ML_Trivial} For a Danielewski surface $\pi:S\rightarrow\mathbb{A}_{k}^{1}$, the following are equivalent :

1) $S$ admits two $\mathbb{A}^{1}$-fibrations with distinct general
fibers.

2) $S$ is isomorphic to a Danielewski surface $S\left(\gamma\right)$
defined by a fine $k$-weighted comb $\gamma=\left(\Gamma,w\right)$.

3) There exists an integer $h\geq1$ and a collection of monic polynomials
$P_{0},\ldots,P_{h-1}\in k\left[t\right]$ with simple roots $a_{i,j}\in k^{*}$,
$i=0,\ldots,h-1$, $j=1,\ldots,\deg_{t}\left(P_{i}\right)$, such
that $S$ is isomorphic to the surface $S_{P_{0},\ldots,P_{h-1}}\subset\textrm{Spec}\left(k\left[x\right]\left[y_{-1},\ldots,y_{h-2}\right]\left[z\right]\right)$
defined by the equations \[
\left\{ \begin{array}{lll}
xz-y_{h-2}{\displaystyle \prod_{l=0}^{h-1}}P_{l}\left(y_{l-1}\right)=0\\
zy_{i-1}-y_{i}y_{h-2}{\displaystyle \prod_{l=i+1}^{h-1}}P_{l}\left(y_{l-1}\right)=0 & xy_{i}-y_{i-1}{\displaystyle \prod_{l=0}^{i}}P_{l}\left(y_{l-1}\right)=0 & 0\leq i\leq h-2\\
y_{i-1}y_{j}-y_{i}y_{j-1}{\displaystyle \prod_{l=i+1}^{j}}P_{l}\left(y_{l-1}\right)=0 &  & 0\leq i<j\leq h-2\end{array}\right.\]

\end{thm}
\begin{proof}
One checks in a similar way as in the proof of Theorem 2.9 in \cite{DubEmb05}
that a surface $S=S_{P_{0},\ldots,P_{h-1}}$ is  a Danielewski surface $\pi={\rm pr}_{x}\mid_{S}:S\rightarrow\mathbb{A}_{k}^{1}$. Furthermore, the projection $\pi'={\rm pr}_{z}\mid_{S}:S\rightarrow\mathbb{A}_{k}^{1}$ is a second $\mathbb{A}^{1}$-fibration on $S$ restricting to
a trivial $\mathbb{A}^{1}$-bundle $\left(\pi'\right)^{-1}\left(\mathbb{A}_{k}^{1}\setminus\left\{ 0\right\} \right)\simeq\textrm{Spec}\left(k\left[z,z^{-1}\right]\left[y_{h-2}\right]\right)$
over $\mathbb{A}_{k}^{1}\setminus\left\{ 0\right\} $. So 3) implies
1). To show that 1) implies 2) we use the following fact, which is a consequence of a result due to M.H. Gizatullin \cite{Giz71} :  if a nonsingular affine surface $S$ defined over an algebraically closed field $k$ admits an $\mathbb{A}^1$-fibration $q:S\rightarrow \mathbb{A}^1_k$, then this fibration is unique up to automorphisms of the base if and only if $S$ does not admit a completion by a nonsingular projective surface $\bar{S}$ for which the boundary divisor  $\bar{S}\setminus S$ is a \emph{zigzag}, that is, a chain of nonsingular proper rational curves.  In \cite{DubG03}, the first author constructed canonical
completions $\bar{S}$ of a Danielewski surface $S\left(\gamma\right)$
defined by a fine $k$-weighted tree $\gamma=\left(\Gamma,w\right)$
for which the dual graph $\Gamma'$ of the boundary divisor $\bar{S}\setminus S\left(\gamma\right)$
is isomorphic to the tree obtained from $\Gamma$ be deleting its
leaves and replacing its root by a chain with two elements. Clearly, $\bar{S}\setminus S\left(\gamma\right)$ is a zigzag  if and only if $\Gamma$ is a comb. The construction given in \emph{loc}. \emph{cit}.
only depends on the existence of an $\mathbb{A}^{1}$-bundle
structure $\rho:S\left(\gamma\right)\rightarrow X\left(r\right)$
on a Danielewski surface $S\left(\gamma\right)$. So it remains valid over
an arbitrary base field $k$. Now let $S=S\left(\gamma\right)$ be a Danielewski surface
admitting two distinct $\mathbb{A}^{1}$-fibrations. Given an algebraic closure $\bar{k}$ of $k$,  the surface $S_{\bar{k}}=S\times_{\textrm{Spec}\left(k\right)}\textrm{Spec}\left(\bar{k}\right)$
is a Danielewski surface isomorphic to the one defined by the
tree $\gamma$ consider as a fine $\bar{k}$-weighted tree via the
inclusion $k\subset\bar{k}$. Since every $\mathbb{A}^{1}$-fibration
$\pi:S\rightarrow\mathbb{A}_{k}^{1}$ lifts to an $\mathbb{A}^{1}$-fibration
$\pi_{\bar{k}}:S_{\bar{k}}\rightarrow\mathbb{A}_{\bar{k}}^{1}$ it
follows that $S_{\bar{k}}$ admits two $\mathbb{A}^{1}$-fibrations
with distinct general fibers. So we deduce from Gizatullin's criterion above
that $\gamma$ is a comb. Thus  1) implies 2).

It remains to show that every Danielewski surface $\pi:S=S\left(\gamma\right)\rightarrow \mathbb{A}^1_k$ defined by a fine $k$-weighted comb $\gamma$ of height $h\geq 1$ admits a closed embedding in an affine space as a surface $S_{P_0,\ldots, P_{h-1}}$. This follows from a general construction described in §4.6 of \cite{DubEmb05} that can be simplified in our more restrictive context. For the convenience of the reader, we indicate below the main steps of the proof. If  $\gamma$ is a chain, then $S\left(\gamma\right)$ is isomorphic to the affine  plane $\mathbb{A}^2_k$ which embeds in $\mathbb{A}^{h+2}_k$ as a surface $S_{P_0,\ldots ,P_{h-1}}$  for which all the polynomials $P_i$, $i=0,\ldots, h-1$ have degree one. We assume from now on that  $\gamma$  has at least two elements at level $1$ (see Remark \ref{rem:Psi_2_Level1} above). We denote by $e_{0,0}<e_{1,0}<\cdots<e_{h-1,0}$ the
elements of the sub-chain $C=\Gamma\setminus L\left(\Gamma\right)$
of $\Gamma$ consisting of elements of $\Gamma$ which are not leaves
of $\Gamma$. For every $l=1,\ldots,h$, the elements of $\Gamma$
at level $l$ distinct from $e_{l,0}$ are denoted by $e_{l,1},\ldots,e_{l,r_{l}}$
provided that they exist. Since $\gamma$ is a comb, it follows from \ref{txt:Abstract_DS_morph} above that $S$ is isomorphic to the surface associated with a certain fine $k$-weighted tree with the same underlying tree $\Gamma$ as $\gamma$ and equipped with a weight function $w$ such that $w\left(\overleftarrow{e_{i,0}e_{i+1,0}}\right)=0$ for every index $i=0,\ldots,h-2$  and such that  $w\left(\overleftarrow{e_{h-1,0}e_{h-1,1}}\right)=0$. We consider $S$ as an  $\mathbb{A}^{1}$-bundle $\rho:S\rightarrow X\left(r\right)$
and we denote by $S_{i}=\textrm{Spec}\left(k\left[x\right]\left[u_{i}\right]\right)$
the trivialising open subsets of $S$ over $X\left(r\right)$. For every
$l=0,\ldots,h-1$ and every $i=1,\ldots,s_{l}$, we let $\tau_{l,i}\left(x,u_{i}\right)=xu_{i}+w\left(\overleftarrow{e_{l-1,0}e_{l,i}}\right)\in k\left[x\right]\left[u_{i}\right]$.
With this notation, the canonical function $\psi$ on $S$
restricts on an open subset $S_{i}$ corresponding to a leaf $e_{l,i}$
of $\Gamma$ at level $l$ to the polynomial $x^{l-1}\tau_{l,i}\left(x,u_{i}\right)\in k\left[x\right]\left[u_{i}\right]$. Therefore, $y_{-1}=\psi$ is constant with
the value $a_{0,i}=w\left(\overleftarrow{e_{0,0}e_{1,i}}\right)\in k^{*}$
on the irreducible component $\pi^{-1}\left(o\right)$ corresponding
to a leaf $e_{1,i}$, $i=1,\ldots,r_{1}$, at level $1$.  It vanishes
identically on every irreducible component of $\pi^{-1}\left(o\right)$
corresponding to a leaf of $\gamma$ at level $l\geq2$. More generally, direct computations show that there exists a unique datum consisting of regular functions $y_{-1},\ldots,y_{h-2}$ and $y_{h-1}$ on $S$ and polynomials $P_i\in k\left[t\right]$, $i=0,\ldots, h-1$ satisfying the following conditions :

a) For every $l=0,\ldots,h-1$, and every $l\leq m\leq h$ , $y_{l-1}$
restricts on an open subset $S_{i}$ corresponding to a leaf $e_{m,i}$
of $\gamma$ at level $m$ to a polynomial $y_{l-1,i}\in k\left[x\right]\left[u_{i}\right]$
such that \begin{eqnarray*}
y_{l-1,i} & = & \left\{ \begin{array}{lll}
L_{l,i}\left(u_{i}\right) & \textrm{mod }x & \textrm{if }m=l\\
a_{l,i}+xL_{l+1,i}\left(u_{i}\right) & \textrm{mod }x^{2} & \textrm{if }m=l+1\\
\xi_{m}x^{m-l-1}\tau_{m,i}\left(x,u_{i}\right)+\nu_{m,i}x^{m-l} & \textrm{mod }x^{m-l+1} & \textrm{if }m>l+1,\end{array}\right.\end{eqnarray*}
where $L_{l,i}\left(u_{i}\right),L_{l+1,i}\left(u_{i}\right)\in k\left[u_{i}\right]$
are  polynomials of degree $1$, $a_{l,i},\xi_{m} \in k^{*}$  and  $\nu_{m,i}\in k$. Furthermore $a_{l,i}\neq a_{l,j}$ for every pair of distinct indices $i$ and $j$.

b) For every $l=0,\ldots,h-1$, $P_{l}$ is the unique monic polynomial
with simple roots $a_{l,1},\ldots,a_{l,r_{l}}$ such that $x^{-1}y_{l-1}\prod_{i=0}^{l-1}P_{i}\left(y_{i-1}\right)P_{l}\left(y_{l-1}\right)$ is a regular function on $S$.

By construction, these functions  $y_{-1},\ldots,y_{h-2},y_{h-1}=z$  distinguish the irreducible
components of the fiber $\pi^{-1}\left(o\right)$ and induce coordinate
functions on them. It follows that the morphism $i=\left(\pi,y_{-1},\ldots,y_{h-1},z\right):S\hookrightarrow\mathbb{A}_{k}^{h+2}$
is an embedding. The same argument as in the proof of Lemma 3.6 in
\cite{DubEmb05} shows that $i$ is actually a closed embedding whose
image is contained in the surface $S_{P_{0},\ldots,P_{h-1}}\subset\mathbb{A}_{k}^{h+2}$
defined in Theorem \ref{thm:Comb_ML_Trivial} above. One checks that the induced morphism $\phi:S\rightarrow S_{P_{0},\ldots,P_{h-1}}$
defines a bijection between the sets of closed points of $S$
and $S_{P_{0},\ldots,P_{h-1}}$. Furthermore, $\phi$ is also birational as $y_{-1}$ induces an isomorphism $\pi^{-1}\left(\mathbb{A}^{1}\setminus\left\{ o\right\} \right)\stackrel{\sim}{\rightarrow}{\rm Spec}\left(k\left[x,x^{-1}\right]\left[y_{-1}\right]\right)$.  Since $S_{P_{0},\ldots,P_{h-1}}$ is nonsingular, we conclude that $\phi$ an isomorphism by virtue of Zariski
Main Theorem (see e.g., 4.4.9 in \cite{EGAIII}).
\end{proof}

\subsection{Special Danielewski surfaces }

\indent\newline\noindent It follows from Adjunction
Formula that every Danielewski surface $S$ in $\mathbb{A}_{k}^{3}$
has a trivial canonical sheaf $\omega_{S/k}=\Lambda^2\Omega^1_{S/k}$. More generally, a Danielewski surface  $\pi:S\rightarrow\mathbb{A}_{k}^{1}$
with a trivial canonical sheaf, or equivalently with
a trivial sheaf of relative differential forms $\Omega_{S/\mathbb{A}_{k}^{1}}^{1}$,
will be called \emph{special}.
\begin{enavant} \label{pro:Spec_DS_charac} These surfaces correspond
to a distinguished class of weighted trees $\gamma$. Indeed, it follows from the gluing construction given in \ref{txt:Abstract_DS_morph} above that a Danielewski surface $S\left(\gamma\right)$ admits a nowhere vanishing differential $2$-form if and only if all the leaves of $\gamma$ are at the same level. In turn,
this means that these surfaces $S$ are the total space of $\mathbb{A}^{1}$-bundles
$\rho:S\rightarrow X\left(r\right)$ over $X\left(r\right)$ defined
by means of transition isomorphisms \[
\tau_{ij}:k\left[x,x^{-1}\right]\left[u_{i}\right]\rightarrow k\left[x,x^{-1}\right]\left[u_{j}\right],\quad u_{i}\mapsto u_{j}+g_{ij}\left(x\right),\quad i,j=1,\ldots,r,\]
 where $g=\left\{ g_{ij}\right\} _{i,j}\in C^1\left(X\left(r\right),\mathcal{O}_{X\left(r\right)}\right)\simeq \mathbb{C}\left[x,x^{-1}\right]^{2r}$ is a \v{C}ech cocycle with
values in the sheaf $\mathcal{O}_{X\left(r\right)}$ for the canonical open
covering $\mathcal{U}_{r}$. So they can be equivalently characterised among Danielewski surfaces by the fact that the underlying $\mathbb{A}^1$-bundle $\rho:S\rightarrow X\left(r\right)$ is actually the structural morphism of a principal homogeneous $\mathbb{G}_a$-bundle.
\end{enavant}
\begin{enavant}To determine isomorphism classes of special
Danielewski surfaces, we can exploit the fact that
the group $\textrm{Aut}\left(X\left(r\right)\right)\simeq\textrm{Aut}\left(\mathbb{A}_{k}^{1}\setminus\left\{ o\right\} \right)\times\mathfrak{S}_{r}$
acts on the set $\mathbb{P}H^{1}\left(X\left(r\right),\mathcal{O}_{X\left(r\right)}\right)$
of isomorphism classes of $\mathbb{A}^{1}$-bundles as above. Indeed, for every $\phi\in\textrm{Aut}\left(X\left(r\right)\right)$,
the image $\phi\cdot\left[g\right]$ of a class $\left[g\right]\in\mathbb{P}H^{1}\left(X\left(r\right),\mathcal{O}_{X\left(r\right)}\right)$
represented by a bundle $\rho:S\rightarrow X\left(r\right)$ is the
isomorphism class of the fiber product bundle ${\rm pr}_{2}:\phi^{*}S=S\times_{X\left(r\right)}X\left(r\right)\rightarrow X\left(r\right)$.
The following criterion generalises a result of J. Wilkens \cite{Wil98}.
\end{enavant}

\begin{thm}
\label{thm:Iso_classes} Two special Danielewski surfaces $\pi_{1}:S_{1}\rightarrow\mathbb{A}_{k}^{1}$
and $\pi_{2}:S_{2}\rightarrow\mathbb{A}_{k}^{1}$ with underlying
$\mathbb{A}^{1}$-bundle structures $\rho_{1}:S_{1}\rightarrow X\left(r_{1}\right)$
and $\rho_{2}:S_{2}\rightarrow X\left(r_{2}\right)$ are isomorphic
as abstract surfaces if and only if $r_{1}=r_{2}=r$ and their isomorphism
classes in $\mathbb{P}H^{1}\left(X\left(r\right),\mathcal{O}_{X\left(r\right)}\right)$
belongs to the same orbit under the action of $\textrm{Aut}\left(X\left(r\right)\right)$.
\end{thm}
\begin{proof}
The condition guarantees that $S_{1}$ and $S_{2}$ are isomorphic.
Suppose conversely that there exists an isomorphism $\Phi:S_{1}\stackrel{\sim}{\rightarrow}S_{2}$.
The divisor class group of a special Danielewski surface $\pi:S\rightarrow\mathbb{A}_{k}^{1}$
is generated by the classes of the connected components $C_{1},\ldots,C_{r}$
of $\pi^{-1}\left(o\right)$ modulo the relation $C_{1}+\cdots+C_{r}=\pi^{-1}\left(o\right)\sim0$,
whence is isomorphic to $\mathbb{Z}^{r-1}$. Therefore, $r_{1}=r_{2}=r$
for a certain $r\geq1$. If one of the $S_{i}$'s, say $S_{1}$ is
isomorphic to a surface $S_{P,1}\subset\mathbb{A}_{k}^{3}$ defined by the equation
$xz-P\left(y\right)=0$, then the result follows from \cite{ML01}.
Otherwise, we deduce from Theorem \ref{thm:Comb_ML_Trivial} that
the $\mathbb{A}^{1}$-fibrations $\pi_{1}:S_{1}\rightarrow\mathbb{A}_{k}^{1}$
and $\pi_{2}:S_{2}\rightarrow\mathbb{A}_{k}^{1}$ are unique up to
automorphisms of the base. In turn, this implies that $\Phi$ induces
an isomorphism $\phi:X\left(r\right)\stackrel{\sim}{\rightarrow}X\left(r\right)$
such that $\phi\circ\rho_{1}=\rho_{2}\circ\Phi$. Therefore, $\Phi:S_{1}\stackrel{\sim}{\rightarrow}S_{2}$
factors through an isomorphism of $\mathbb{A}^{1}$-bundles $\tilde{\phi}:S_{1}\stackrel{\sim}{\rightarrow}\phi^{*}S_{2}$,
where $\phi^{*}S_{2}$ denotes the the fiber product $\mathbb{A}^{1}$-bundle
${\rm pr}_{2}:\phi^{*}S_{2}=S_{2}\times_{X\left(r\right)}X\left(r\right)\rightarrow X\left(r\right)$.
This completes the proof as $\phi^{*}S_{2}\simeq S_{2}$.
\end{proof}

\section{Danielewski surfaces in $\mathbb{A}_{k}^{3}$ defined by an equation
of the form $x^{h}z-Q\left(x,y\right)=0$ and their automorphisms}

In this section, we study Danielewski surfaces $\pi:S\rightarrow\mathbb{A}_{k}^{1}$ non isomorphic to $\mathbb{A}^2_k$  admitting a closed embedding $i:S\hookrightarrow\mathbb{A}_{k}^{3}$ in the affine $3$-space as a surface $S_{Q,h}$ defined by the equation $x^{h}z-Q\left(x,y\right)=0$. We show that a same abstract Danielewski surface may admit many such closed embeddings. In particular, we establish that $S$ can be embedded as a surface $S_{\sigma,h}$ defined by an equation of the form  $x^{h}z-\prod_{i=1}^r\left(y-\sigma_i\left(x\right)\right)=0$ for a suitable collection of polynomials $\sigma=\left\{\sigma_i\left(x\right)\right\}_{i=1,\ldots,r}$. Next we study the automorphism groups of the above surfaces $S$. We show that, in a closed embedding as a surface $S_{\sigma,h}$, every automorphism of $S$ explicitly arises as the restriction of an automorphism of the ambient space.
We will show on the contrary in the next section that it is not true for a general embedding as a surface $S_{Q,h}$.

\subsection{Danielewski surfaces $S_{Q,h}$}
\indent\newline\noindent  A surface $S=S_{Q,h}$ in $\mathbb{A}_{k}^{3}$
defined by the equation $x^{h}z-Q\left(x,y\right)=0$ is a Danielewski surface
$\pi={\rm pr}_{x}\mid_{S}:S\rightarrow\mathbb{A}_{k}^{1}$
if and only if the polynomial $Q\left(0,y\right)$ splits with simple
roots $y_{1},\ldots,y_{r}\in k$, where $r=\deg_{y}\left(Q\left(0,y\right)\right)$.
If $r=1$, then $\pi^{-1}\left(o\right)\simeq\mathbb{A}_{k}^{1}$
and $\pi:S\rightarrow\mathbb{A}_{k}^{1}$ is isomorphic to a trivial
$\mathbb{A}^{1}$-bundle. Thus $S$ is isomorphic to the affine plane.
Otherwise, if $r\geq2$, then $S$ is not isomorphic to $\mathbb{A}_{k}^{2}$,
as follows for instance from the fact that the divisor class group
$\textrm{Div}\left(S\right)$ of $S$ is isomorphic to $\mathbb{Z}^{r-1}$,
generated by the classes of the connected components $C_{1},\ldots,C_{r}$
of $\pi^{-1}\left(o\right)$, with a unique relation $C_{1}+\ldots+C_{r}=\textrm{div}\left(\pi^{*}x\right)\sim0$.

 The above class of Danielewski surfaces contains affine surfaces  $S_{P,h}$ in $\mathbb{A}^3_k$ defined by an equation of the form $x^hz-P\left(y\right)=0$, where $P\left(y\right)$ is a polynomial which splits with simple roots $y_1,\ldots, y_r$ in $k$. Replacing the constants $y_{i}\in k$ by suitable polynomials  $\sigma_{i}\left(x\right)\in k\left[x\right]$ leads to the following more general class of examples.

\begin{example}
\label{exa:Main_example} Let $h\geq1$  be an integer and let $\sigma=\left\{ \sigma_{i}\left(x\right)\right\} _{i=1,\ldots,r}$
be a collection of $r\geq2$ polynomials $\sigma_{i}\left(x\right)=\sum_{j=0}^{h-1}\sigma_{i,j}x^{j}\in k\left[x\right]$ such that
$\sigma_{i}\left(0\right)\neq\sigma_{j}\left(0\right)$ for every
$i\neq j$. The surface $S=S_{\sigma,h}$ in $\mathbb{A}_{k}^{3}={\rm Spec}\left(k\left[x,y,z\right] \right)$
defined by the equation \[x^{h}z-\prod_{i=1}^{r}\left(y-\sigma_{i}\left(x\right)\right)=0\]
is a Danielewski surface $\pi={\rm pr}_{x}\mid_{S}:S\rightarrow\mathbb{A}_{k}^{1}$.
The fiber $\pi^{-1}\left(o\right)$ consists of $r$ copies
$C_{i}$ of the affine line defined by the equations $\left\{ x=0,y=\sigma_{i}\left(0\right)\right\} _{i=1,\ldots,r}$ respectively.
For every index $i=1,\ldots,r$, the open subset $S_{i}=S\setminus\bigcup_{j\neq i}C_{i}$
of $S$ is isomorphic to the affine plane $\mathbb{A}_{k}^{2}=\textrm{Spec}\left(k\left[x,u_{i}\right]\right)$,
where $u_{i}$ denotes the regular function on $S_{i}$ induced by
the rational function \[ u_{i}=x^{-h}\left(y-\sigma_{i}\left(x\right)\right)=z\prod_{j\neq i}\left(y-\sigma_{j}\left(x\right)\right)^{-1}\in k\left(S\right) \]
on $S$. It follows that $\pi:S\rightarrow\mathbb{A}_{k}^{1}$
factors through an $\mathbb{A}^{1}$-bundle $\rho:S\rightarrow X\left(r\right)$
isomorphic to the one with transition pairs $\left(f_{ij},g_{ij}\right)=\left(1,x^{-h}\left(\sigma_{j}\left(x\right)-\sigma_{i}\left(x\right)\right)\right)$,
$i,j=1,\ldots,r$.
The collection $\sigma=\left\{\sigma_i\left(x\right)\right\}_{i=1,\ldots,r}$ is exactly the one associated with the following fine $k$-weighted tree $\gamma=\left(\Gamma,w\right)$.

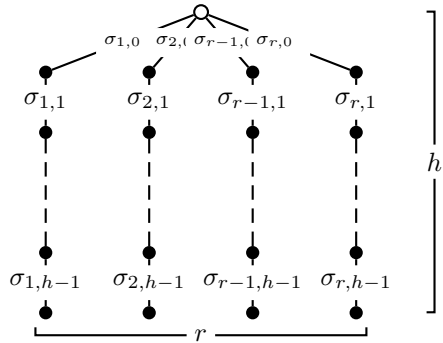
\begin{pspicture}(-4.6,2.5)(8,-2.7)

\def\dedge{\ncline[linestyle=dashed]}

\rput(2,0){

\pstree[treemode=D,radius=2.5pt,treesep=1.2cm,levelsep=0.8cm]{\Tc{3pt}}{

  \pstree{\TC*\mput*{{\scriptsize $\sigma_{1,0}$}}} {

              \pstree{\TC*\mput*{$\sigma_{1,1}$}}{ \skiplevels{1}

                 \pstree{\TC*[edge=\dedge]}{

                    \TC*\mput*{$\sigma_{1,h-1}$}

                                                               }

                                                      \endskiplevels

                                                     }

                                          }

\pstree{\TC*\mput*{{\scriptsize $\sigma_{2,0}$}}}{

            \pstree{\TC*\mput*{$\sigma_{2,1}$}}{\skiplevels{1}

                     \pstree{\TC*[edge=\dedge]} {

                              \TC*\mput*{$\sigma_{2,h-1}$}

                                                                    }

                                                          \endskiplevels

                                                           }

                                                }

\pstree{\TC*\mput*{{\scriptsize $\sigma_{r-1,0}\;$}}} {

              \pstree{\TC*\mput*{$\sigma_{r-1,1}$}}{ \skiplevels{1}

                 \pstree{\TC*[edge=\dedge]}{

                    \TC*\mput*{$\sigma_{r-1,h-1}$}

                                                               }

                                                      \endskiplevels

                                                     }

                                          }

\pstree{\TC*\mput*{{\scriptsize $\sigma_{r,0}$} }}{

            \pstree{\TC*\mput*{$\sigma_{r,1}$}}{\skiplevels{1}

                     \pstree{\TC*[edge=\dedge]} {

                              \TC*\mput*{$\sigma_{r,h-1}$}

                                                                    }

                                                          \endskiplevels

                                                           }

                                                }

}

}

\pnode(-0.2,-2.2){A}\pnode(4.2,-2.2){B}

\ncbar[angleA=270, arm=3pt]{A}{B}\ncput*[npos=1.5]{$r$}

\pnode(5,2){C}\pnode(5,-2){D}

\ncbar[arm=3pt]{C}{D}\ncput*[npos=1.5]{$h$}

\end{pspicture}

\noindent So $S$ is isomorphic to the corresponding Danielewski surface $\pi_{\gamma}:S\left(\gamma\right)\rightarrow \mathbb{A}^1_k$. By definition (see \ref{txt:canonical_morphism} above),  the canonical function $\psi_{\gamma}$ on $S\left(\gamma\right)$ is the unique regular function restricting to the polynomial function $\psi_{\gamma,i}=x^hu_i+\sigma_i\left(x\right)\in k\left[x,u_i\right]$ on the trivialising open subsets $S_i\simeq\mathbb{A}^2_k$, $i=1,\ldots,r$ of $S\left(\gamma\right)$. So it coincides with the restriction of $y$ on $S$ under the above isomorphism. In the setting of Theorem \ref{thm:GenDanMor_2_Tree}, this means that  $\gamma$ corresponds to the Danielewski surface $S$ equipped with the birational morphism ${\rm pr}_{x,y}:S\rightarrow \mathbb{A}^2_k$.
\end{example}

It turns out that up to isomorphisms, the above class of Danielewski surfaces $S_{\sigma,h}$ contains all possible Danielewski surfaces  $S_{Q,h}$, as shown by the following result.

\begin{thm}
\label{thm:Equivalent_charac} Let $S_{Q,h}$ be a  Danielewski surface  in $\mathbb{A}^3_k$ defined by the equation $x^hz-Q\left(x,y\right)=0$, where $Q\left(x,y\right)\in k\left[x,y\right]$ is a polynomial such that $Q\left(0,y\right)$ splits with $r\geq 2$ simples roots in $k$. Then there exists a  collection   $\sigma=\left\{\sigma_i\left(x\right)\right\}_{i=1,\ldots, r}$ of polynomials of degrees $\deg\left(\sigma_i\left(x\right)\right)< h$ such that $S_{Q,h}$ is isomorphic to the surface $S_{\sigma,h}$ defined by the equation
$x^hz-\prod_{i=1}^r\left(y-\sigma_i\left(x\right)\right)=0$.
\end{thm}

\begin{proof}
Since $Q\left(0,y\right)$ splits with simple roots  $y_1,\ldots, y_r$ in $k$, a variant of the classical Hensel Lemma (see e.g., Theorem 7.18 p. 208 in \cite{Eis95}) guarantees that the polynomial $Q(x,y)$ can be written in a unique way as \[Q\left(x,y\right)=R_{1}\left(x,y\right)\prod_{i=1}^{r}\left(y-\sigma_{i}\left(x\right)\right)+x^{h}R_{2}\left(x,y\right), \]
where  $R_{1}\left(x,y\right)\in k\left[x,y\right]\setminus\left(x^{h}k\left[x,y\right]\right)$
is a polynomial such that $R_{1}\left(0,y\right)$ is a nonzero constant and where $\sigma=\left\{\sigma_i\left(x\right)\right\}_{i=1,\ldots, r}$ is a collection of polynomials of degree strictly lower than $h$ such that $\sigma_i\left(0\right)=y_i$ for every index $i=1,\ldots, r$. Since $y_i\neq y_j$ for every $i\neq j$ and $R_1\left(0,y\right)$ is a nonzero constant, it follows that for every index $i=1,\ldots,r$, the rational function
\[ u_i=x^{-h} \left(y-\sigma_i\left(x\right)\right)=\prod_{j\neq i}\left(y-\sigma_j\left(x\right)\right)^{-1}R_1\left(x,y\right)^{-1}\left(z-R_2\left(x,y\right)\right)\]
on $S_{Q,h}$ restricts to a regular function on the complement $S_i$ in $S_{Q,h}$ of the irreducible components of the fiber ${\rm pr}_x^{-1}\left(0\right)$ defined by the equations $\left\{x=0,y=y_j\right\}_{j\neq i}$ and induces  an isomorphism $S_i\simeq{\rm Spec}\left(k\left[x,u_i\right]\right)$.
Therefore, the collection $\sigma=\left\{\sigma_i\left(x\right)\right\}_{i=1,\ldots,r}$ is precisely the one associated with the fine $k$-weighted rake $\gamma=\left(\Gamma,w\right)$ with all
its leaves at a same level $h$ corresponding to the Danielewski surface ${\rm pr}_x:S_{Q,h}\rightarrow \mathbb{A}^1_k$ equipped with the birational morphism $\psi={\rm pr}_{x,y}:S_{Q,h}\rightarrow\mathbb{A}^2_k$ (see \ref{thm:GenDanMor_2_Tree} and \ref{pro:Spec_DS_charac} above). In turn, we deduce from example \ref{exa:Main_example} that the Danielewski surface $S\left(\gamma\right)$ associated with $\gamma $ embeds as the surface $S_{\sigma,h}$ in $\mathbb{A}^3_k$ defined by the equation $x^hz-\prod_{i=1}^r\left(y-\sigma_i\left(x\right)\right)=0$. This completes the proof.
\end{proof}

\begin{defn}
\label{def:Embed_def}
Given a Danielewski surface  $S$ isomorphic to a certain surface $S_{Q,h}$ in $\mathbb{A}^3_k$, a closed  embedding $i_s:S\hookrightarrow\mathbb{A}_{k}^{3}$  of $S$ in $\mathbb{A}^3_k$  as a surface $S_{\sigma,h}$ defined by the equation \[x^hz-\prod_{i=1}^r\left(y-\sigma_i\left(x\right)\right)=0\] is  called  a \emph{standard embedding of} $S$. We  say that $S_{\sigma,h}$ is a \emph{standard
form} \emph{of} $S$ \emph{in} $\mathbb{A}_{k}^{3}$.
\end{defn}

\begin{enavant} \label{lem:Hensel_lemma} \label{rem:def phi_s}
It follows from the above discussion that every Danielewski surface $S$ isomorphic to a certain surface $S_{Q,h}$ in $\mathbb{A}^3_k$ admits a standard embedding in $\mathbb{A}^3_k$. Following the proof of Theorem \ref{thm:Equivalent_charac}, we can in fact construct explicitly the isomorphisms between a Danielewski surface $S_{Q,h}$ and one of its standard forms $S_{\sigma,h}$. Let $Q(x,y)=R_{1}(x,y)\prod_{i=1}^{r}\left(y-\sigma_{i}\left(x\right)\right)+x^{h}R_{2}(x,y)$ be as in the proof of Theorem \ref{thm:Equivalent_charac}. Then, the endomorphism $\Phi^s$ of $\mathbb{A}^3_k$ defined by $\left(x,y,z\right)\mapsto\left(x,y,R_{1}\left(x,y\right)z+R_{2}\left(x,y\right)\right)$ induces an isomorphism $\phi^s$ between $S_{\sigma,h}$ and $S_{Q,h}$. One checks conversely that for every pair $\left(f,g\right)$ of polynomials such that
$R_{1}\left(x,y\right)f\left(x,y\right)+x^{h}g\left(x,y\right)=1$, the endomorphism $\Phi_s$ of $\mathbb{A}^3_k$ defined by \[\left(x,y,z\right)\mapsto\left(x,y,f\left(x,y\right)z+g\left(x,y\right)\prod_{i=1}^{r}\left(y-\sigma_{i}\left(x\right)\right)-f\left(x,y\right)R_{2}\left(x,y\right)\right)\] induces an isomorphism $\phi_s$ between $S_{Q,h}$ and $S_{\sigma,h}$ such that $\phi^s\circ\phi_s={\rm id}_{S_{Q,h}}$ and $\phi_s\circ\phi^s={\rm id}_{S_{\sigma,h}}$. Note that  since $R_1\left(0,y\right)$ is a nonzero constant, the regular function $\xi=x^{-h}(R_{1}\prod_{i=1}^{r}\left(y-\sigma_{i}\left(x\right)\right))+R_{2}$ on $S_{\sigma,h}$ still induces a coordinate function on every irreducible component of the fiber $\pi^{-1}\left(o \right)$ of the morphism $\pi={\rm pr}_x:S_{\sigma,h}\rightarrow \mathbb{A}^1_k$, and  the regular functions $\pi$, $y$ and $\xi$  define a new closed embedding of $S_{\sigma,h}$ in $\mathbb{A}_{k}^{3}$ inducing an isomorphism between $S_{\sigma,h}$ and the surface $S_{Q,h}$. This can be interpreted by saying that a closed embedding $i_{Q,h}:S\hookrightarrow\mathbb{A}_{k}^{3}$ of a Danielewski surface $S$ in $\mathbb{A}_{k}^{3}$ as a surface $S_{Q,h}$ is a twisted form of a standard embedding of $S$ obtained by modifying the function inducing a coordinate on every irreducible component of the fiber $\pi^{-1}\left(o\right)$. 
\end {enavant}

\begin{enavant} Using standard forms makes the study of isomorphism classes of Danielewski surfaces $S_{Q,h}$ simpler. For instance, we have the following characterisation which generalises a result due to L. Makar-Limanov
\cite{ML01} for complex surfaces $S_{P,h}$ defined by the equations $x^{h}z-P\left(y\right)=0$.
\end{enavant}

\begin{prop} \label{thm:Normal_forms_iso} Two Danielewski surfaces $S_{\sigma_1,h_1}$ and $S_{\sigma_2,h_2}$  in $\mathbb{A}_{k}^{3}$ defined by the  equations \[
x^{h_{1}}z=P_{1}\left(x,y\right)=\prod_{i=1}^{r_{1}}\left(y-\sigma_{1,i}\left(x\right)\right)\quad\textrm{and}\quad x^{h_{2}}z=P_{2}\left(x,y\right)=\prod_{i=1}^{r_{2}}\left(y-\sigma_{2,i}\left(x\right)\right)\]
 are isomorphic if and only if $h_{1}=h_{2}=h$, $r_{1}=r_{2}=r$
 and there exists a triple $\left(a,\mu,\tau\left(x\right)\right)\in k^{*}\times k^{*}\times k\left[x\right]$
such that $P_{2}\left(ax,y\right)=\mu^{r}P_{1}\left(x,\mu^{-1}y+\tau\left(x\right)\right)$.
\end{prop}
\begin{proof}
The condition is sufficient. Indeed, one checks that the automorphism  \[\left(x,y,z\right)\mapsto \left(ax,\mu \left(y-\tau\left(x\right)\right), \mu^r a^{-2}z\right)\] of $\mathbb{A}^3_k$ induces an isomorphism between $S_{\sigma_1,h}$ and $S_{\sigma_2,h}$. Conversely,
suppose that  $S_1=S_{\sigma_1,h_1}$ and $S_2=S_{\sigma_2,h_2}$  are isomorphic. Then $h_1=h_2=h$ and $r_1=r_2=r$ by virtue of Theorem \ref{thm:Iso_classes} above.  If $h=1$ then the result follows from \cite{ML01}. Otherwise, if $h\geq 2$ then  it follows from Theorem \ref{thm:Comb_ML_Trivial} and example \ref{exa:Main_example} above that the underlying $\mathbb{A}^{1}$-bundle structures $\rho_1:S_1\rightarrow X\left(r\right)$ and $\rho_2:S_2\rightarrow X\left(r\right)$ corresponding to the transition functions
\[\left\{g_{1,ij}=x^{-h}\left(\sigma_{1,j}\left(x\right)-\sigma_{1,i}\left(x\right)\right)\right\}_{i,j=1,\ldots,r} \textrm{ and } \left\{g_{2,ij}=x^{-h}\left(\sigma_{2,j}\left(x\right)-\sigma_{2,i}\left(x\right)\right)\right\}_{i,j=1,\ldots,r}\]
respectively are unique such structures on $S_1$ and $S_2$ up to automorphisms of the base $X\left(r\right)$. Therefore, every isomorphism $\Phi:S_1\stackrel{\sim}{\rightarrow} S_2$ induces an automorphism $\phi$ of $X\left(r\right)$ such that $\rho_2\circ\Phi=\phi\circ \rho_1$. Consequently, every such isomorphism $\Phi$  is  determined by a collection of local isomorphisms $\Phi_{i}:S_{1,i}\stackrel{\sim}{\rightarrow}S_{2,\alpha\left(i\right)}$ where $\alpha\in\mathfrak{S}_{r}$, defined by $k$-algebra isomorphisms
\[
\Phi_{i}^{*}:k\left[x\right]\left[u_{2,\alpha\left(i\right)}\right]\longrightarrow k\left[x\right]\left[u_{1,i}\right],\quad x\mapsto a_{i}x,\quad u_{2,\alpha\left(i\right)}\mapsto\lambda_{i}u_{1,i}+b_{i}\left(x\right),\quad i=1,\ldots,r\]
 where $a_{i},\lambda_{i}\in k^{*}$ and where $b_{i}\in k\left[x\right]$. These local isomorphisms glue to a global one if and only if  $a_{i}=a$ and $\lambda_{i}=\lambda$ for every index $i=1,\ldots , r$,
and the relation $\lambda g_{1,ij}\left(x\right)+b_{i}\left(x\right)=g_{2,\alpha\left(i\right)\alpha\left(j\right)}\left(ax\right)+b_{j}\left(x\right)$
holds in $k\left[x,x^{-1}\right]$ for every indices $i,j=1,\ldots,r$. Since the $\sigma_{1,i}$'s and $\sigma_{2,i}$'s have degrees strictly lower than $h$, we conclude that the latter condition is equivalent to the fact that $b_{i}\left(x\right)=b\left(x\right)$ for every $i=1,\ldots,r$ and that
the polynomial $c\left(x\right)=\sigma_{2,\alpha\left(i\right)}\left(ax\right)-\lambda a^{h}\sigma_{1,i}\left(x\right)$ does not depend on the index $i$. Letting $\mu=\lambda a^h$ and $\tau\left(x\right)=\mu^{-1}c\left(x\right)$, this means exactly  that $P_2\left(ax,y\right)=\mu^{r}P_1\left(x,\mu^{-1}y+\tau\left(x\right)\right)$.
\end{proof}

\begin{enavant}  \label{rem:standard embeddings} \label{txt:Non_algeb_equiv}  The proof above implies in particular that all  standard embeddings of a same Danielewski surface are algebraically equivalent. It is natural to ask if a closed embedding $i_{Q,h}:S\hookrightarrow \mathbb{A}^3_k$ of Danielewski surface $S$ as a surface $S_{Q,h}$ is algebraically equivalent to a standard one.If so, then we say that the embedding $i_{Q,h}$ is \emph{rectifiable}. The fact that the endomorphisms $\Phi^s$ and $\Phi_s$ of $\mathbb{A}^3_k$ constructed in \ref{rem:def phi_s} are not invertible in general may lead one to suspect that there exists non-rectifiable embeddings of Danielewski surfaces nonisomorphic to the affine plane. This is actually the case, and the first known examples have been  recently discovered by  G. Freudenburg and L. Moser-Jauslin \cite{FrMo02}. For instance, they  established that the surface $S_1$ in $\mathbb{A}^3_{\mathbb{C}}$ defined by the equation $f_1=x^2z-\left(1-x\right)\left(y^2-1\right)=0$ is a non-rectifiable  embedding of a Danielewski surface. Indeed, a standard form for $S_1$ would be the Danielewski surface $S_0$ defined by the equation $f_0=x^2z-\left(y^2-1\right)=0$.  We observe that the level surface $f_{0}^{-1}\left(1\right)$ of $f_{0}$ is a singular
surface. On the other hand, all the level surfaces of $f_{1}$ are
nonsingular as follows for instance from the Jacobian Criterion. Therefore, condition 3) in Definition \ref{def:Algebraic_equiv_def} cannot be satisfied and so, it is impossible to find an automorphism of $\mathbb{A}^3_{\mathbb{C}}$ mapping $S_1$ isomorphically onto $S_0$.

The classification of these embeddings up to algebraic equivalence is a difficult problem in general (see \cite{MoP05} for the case $h=r=2$).
However, if $k=\mathbb{C}$, the following result shows that things become simpler if one works in the holomorphic category.
\end{enavant}

\begin{thm}
\label{thm:Embed_analytic_equiv} The embeddings $i_{Q,h}:S\hookrightarrow\mathbb{A}_{\mathbb{C}}^{3}$
of a Danielewski surface $S$ as a surface defined by the equation $x^{h}z-Q\left(x,y\right)=0$
are all \emph{analytically} equivalent.
\end{thm}
\begin{proof}
It suffices to show that every embedding $i_{Q,h}$ is analytically
equivalent to a standard one $i_{\sigma,h}$. In view of the proof of Theorem \ref{thm:Equivalent_charac}, we can let  $Q\left(x,y\right)=R_1\left(x,y\right)\prod_{i=1}^{r}\left(y-\sigma_{i}\left(x\right)\right)+x^hR_2(x,y)$.
It is enough  to construct an holomorphic automorphism $\Psi$ of $\mathbb{A}_{\mathbb{C}}^{3}$
such that $$\Psi^{*}\left(x^{h}z-\prod_{i=1}^{r}\left(y-\sigma_{i}\left(x\right)\right)\right)=\alpha\left(x^{h}z-Q\left(x,y\right)\right)$$
for a suitable invertible holomorphic function $\alpha$ on $\mathbb{A}_{\mathbb{C}}^{3}$.
We let $R_1\left(0,y\right)=\lambda\in\mathbb{C}^{*}$ and we let $f\left(x,y\right)\in\mathbb{C}\left[x,y\right]$
be a polynomial such that $\lambda\exp\left(xf\left(x,y\right)\right)\equiv R_1\left(x,y\right)$
mod $x^{h}$. Now the result follows from the fact that the holomorphic
automorphism $\Psi$ of $\mathbb{A}_{\mathbb{C}}^{3}$ defined by

$$\Psi\left(x,y,z\right)= \left(x,y,\lambda\exp\left(xf\left(x,y\right)\right)z-x^{-h}[\lambda\exp\left(xf\left(x,y\right)\right)-R_1\left(x,y\right)]\prod_{i=1}^{r}\left(y-\sigma_{i}\left(x\right)\right)+R_2(x,y)\right)$$
 satisfies $\Psi^{*}\left(x^{h}z-Q\left(x,y\right)\right)=\lambda\exp\left(xf\left(x,y\right)\right)\left(x^{h}z-\prod_{i=1}^{r}\left(y-\sigma_{i}\left(x\right)\right)\right)$.
\end{proof}

\begin{example}
We observed in \ref{txt:Non_algeb_equiv} that the
surfaces $S_{0}$ and $S_{1}$ defined by the equations $f_0=x^2z-\left(y^2-1\right)=0$
and $f_1=x^2z-\left(1-x\right)\left(y^2-1\right)=0$ are algebraically
inequivalent embeddings of a same surface $S$. However, they are
analytically equivalent via the automorphism $\left(x,y,z\right)\mapsto\left(x,y,e^{-x}z-x^{-2}\left(e^{-x}-1+x\right)(y^2-1)\right)$
of $\mathbb{A}_{\mathbb{C}}^{3}$.
\end{example}

\subsection{Automorphisms of Danielewski surfaces $S_{Q,h}$ in $\mathbb{A}_{k}^{3}$ }

\indent\newline\noindent In \cite{ML90} and \cite{ML01}, Makar-Limanov
computed the automorphism groups of surfaces in $\mathbb{A}^{3}$
defined by the equation $x^{h}z-P\left(y\right)=0$, where $h\geq1$ and where
$P\left(y\right)$ is an arbitrary polynomial. In particular, he established
that every automorphism of such a surface is induced by the restriction
of an automorphism of the ambient space. Recently, A. Crachiola \cite{Cra06} established
that this also holds for surfaces defined by the equations $x^{h}z-y^{2}-r\left(x\right)y=0$,
where $h\geq1$ and where $r\left(x\right)$ is an arbitrary polynomial
such that $r\left(0\right)\neq0$. This subsection is devoted to the proof of the more general structure Theorem \ref{thm:S_Q,h autos} below. We begin with the case of Danielewski surfaces in standard form.


\begin{thm}
\label{thm:Main_auto_thm} The automorphism group of a Danielewski surface
 $S_{\sigma,h}$ defined by the equation \[ x^hz-P\left(x,y\right)=0,\qquad \textrm{ where   } \quad  P\left(x,y\right)=\prod_{i=1}^r\left(y-\sigma_i\left(x\right)\right)\]   is induced by the restriction  of an automorphism of $\mathbb{A}^3_k$  belonging to the subgroup $G_{\sigma,h}$ of ${\rm Aut}\left(\mathbb{A}_{k}^{3}\right)$
generated by the following automorphisms:

\emph{(a)} $\Delta_{b}\left(x,y,z\right)=\left(x,y+x^{h}b\left(x\right),z+x^{-h}\left(P\left(x,y+x^{h}b\left(x\right)\right)-P\left(x,y\right)\right)\right)$,
where $b\left(x\right)\in k\left[x\right]$.

\emph{(b)} If there exists a polynomial $\tau\left(x\right)$ such
that $P\left(x,y+\tau\left(x\right)\right)=\tilde{P}\left(y\right)$
then the automorphisms $H_{a}\left(x,y,z\right)=\left(ax,y+\tau\left(ax\right)-\tau\left(x\right),a^{-h}z\right)$,
where $a\in k^{*}$ should be added.

\emph{(c)} If there exists a polynomial $\tau\left(x\right)$ such
that $P\left(x,y+\tau\left(x\right)\right)=\tilde{P}\left(x^{q_{0}},y\right)$,
then the cyclic automorphisms $\tilde{H}_{a}\left(x,y,z\right)=\left(ax,y+\tau\left(ax\right)-\tau\left(x\right),a^{-h}z\right)$,
where $a\in k^{*}$ and $a^{q_{0}}=1$ should be added.

\emph{(d)} If there exists a polynomial $\tau\left(x\right)$ such
that $P\left(x,y+\tau\left(x\right)\right)=y^{i}\tilde{P}\left(x,y^{s}\right)$,
where $i=0,1$ and $s\geq2$, then the cyclic automorphisms $S_{\mu}\left(x,y,z\right)=\left(x,\mu y+\left(1-\mu\right)\tau\left(x\right),\mu^{i}z\right)$,
where $\mu\in k^{*}$ and $\mu^{s}=1$ should be added.

\emph{(e)} If $\textrm{char}\left(k\right)=s>0$ and $P\left(x,y\right)=\tilde{P}\left(y^{s}-c\left(x\right)^{s-1}y\right)$
for a certain polynomial $c\left(x\right)\in k\left[x\right]$ such
that $c\left(0\right)\neq0$, then the automorphism $T_{c}\left(x,y,z\right)=\left(x,y+c\left(x\right),z\right)$
should be added.

\emph{(}f\emph{)} If $h=1$,  then the involution $I\left(x,y,z\right)=\left(z,y,x\right)$
should be added.
\end{thm}

\begin{rem}\label{rem: k+-actions}
Automorphisms of type a) in Theorem \ref{thm:Main_auto_thm} correspond to algebraic actions of the additive group $\mathbb{G}_a$ on the surface $S_{\sigma,h}$. More precisely, for every polynomial $b\in k\left[x\right]$, the subgroup $\left\{\Delta_{tb\left(x\right)}, t\in k \right\}$ of ${\rm Aut}\left(S_{\sigma,h}\right)$ is isomorphic to $\mathbb{G}_a$, the corresponding $\mathbb{G}_a$-action on $S_{\sigma,h}$ being defined by $t\star\left(x,y,z\right)=\Delta_{tb\left(x\right)}\left(x,y,z\right)$. 
Similarly, automorphisms of type b) correspond to algebraic actions of the multiplicative group $\mathbb{G}_m$. 
\end{rem}

\begin{proof} It is clear that every automorphism of $\mathbb{A}_{k}^{3}$ of types
(a)-(f) above leaves $S_{\sigma,h}$ invariant, whence induces an
automorphism of $S_{\sigma,h}$. If $h=1$, then the converse
follows from \cite{ML90}. Otherwise, if $h\geq2$, then the same argument as the one used in the proof of Proposition \ref{thm:Normal_forms_iso} above show  that every automorphism of $S_{\sigma,h}$
is determined by a datum $\mathcal{A}_{\Phi}=\left(\alpha,\mu,a,b\left(x\right)\right)$
such that that the polynomial $c\left(x\right)=\sigma_{\alpha\left(i\right)}\left(ax\right)-\mu\sigma_{i}\left(x\right)+x^{h}b\left(x\right)$
does not depend on the index $i=1,\ldots,r$. Furthermore, it follows from the construction of the closed embedding  of $S_{\sigma,h}$ in $\mathbb{A}^3_k$ given in Example \ref{exa:Main_example} that every such collection correspond to an automorphism of $S_{\sigma,h}$  induced by the restriction of the following
automorphism $\Psi$ of $\mathbb{A}_{k}^{3}$: \[
\Psi\left(x,y,z\right)=\left(ax,\mu y+c\left(x\right),a^{-h}\mu^rz+\left(ax\right)^{-h}(\prod_{i=1}^{r}\left(\mu y+c\left(x\right)-\sigma_{i}\left(ax\right)\right)-\mu^r\prod_{i=1}^{r}\left(y-\sigma_{i}\left(x\right)\right))\right).\]
One checks easily using this description that the composition of  two automorphisms $\Phi_{1}$ and $\Phi_{2}$ of $S_{\sigma,h}$ defined by data $\mathcal{A}_{\Phi_{1}}=\left(\alpha_{1},\mu_{1},a_{1},b_{1}\right)$ and $\mathcal{A}_{\Phi_{2}}=\left(\alpha_{2},\mu_{2},a_{2},b_{2}\right)$
is the automorphism  with corresponding datum $\mathcal{A}_{\Phi}=\left(\alpha_{2}\circ\alpha_{1},\mu_{2}\mu_{1},a_{2}a_{1},a_2^{-h}\mu_{2}b_{1}\left(x\right)+b_{2}\left(a_{1}x\right)\right)$.

Clearly, automorphisms of type (a) coincide with the
ones determined by  data $\mathcal{A}=\left(\textrm{Id},1,1,b\left(x\right)\right)$, where $b\left(x\right)\in k\left[x\right]$. In view of the composition rule above, it suffices to consider from now on automorphisms corresponding to
data $\mathcal{A}=\left(\alpha,\mu,a,0\right)$.

1°) If $\alpha$ is trivial, then $\mu=1$ by virtue of Lemma
\ref{pro:auto_data} below, and so $\mathcal{A}=\left(\textrm{Id},1,a,0\right)$. Then, the relation $c\left(x\right)=\sigma_{i}\left(ax\right)-\sigma_{i}\left(x\right)$ holds for every $i=1,\ldots,r$.

1°a) If $a^{q}\neq1$ for every $q=1,\ldots,h-1$, then there exists
a polynomial $\tau\left(x\right)\in k\left[x\right]$ such that $\sigma_{i}\left(x\right)=\sigma_{i}\left(0\right)+\tau\left(x\right)$
for every $i=1,\ldots,r$. Thus $c\left(x\right)=\tau\left(ax\right)-\tau\left(x\right)$
and $P\left(x,y+\tau\left(x\right)\right)=\tilde{P}\left(y\right)=\prod_{i=1}^{r}\left(y-\sigma_{i}\left(0\right)\right)$ and the  corresponding automorphism is of type (b).

1°b) If $a\neq1$ but $a^{q_{0}}=1$ for a minimal $q_{0}=2,\ldots,h-1$,
then there exists polynomials $\tau\left(x\right)$ and $\tilde{\sigma}_{i}\left(x\right)$,
$i=1,\ldots,r$, such that $\sigma_{i}\left(x\right)=\tilde{\sigma}_{i}\left(x^{q_{0}}\right)+\tau\left(x\right)$
for every $i=1,\ldots,r$. So there exists a polynomial $\tilde{P}$
such that $P\left(x,y+\tau\left(x\right)\right)=\tilde{P}\left(x^{q_{0}},y\right)$.
Moreover, $c\left(x\right)=\tau\left(ax\right)-\tau\left(x\right)$
and the corresponding automorphism is of type (c).

2°) If $\alpha$ is not trivial then $\mu^{s}=1$. Since $\Phi=\Phi_{2}\circ\Phi_{1}$,
where $\Phi_{1}$and $\Phi_{2}$ denote the automorphisms with data
$\mathcal{A}_{\Phi_{1}}=\left(\textrm{Id},1,a,0\right)$ and $\mathcal{A}_{\Phi_{2}}=\left(\alpha,\mu,1,0\right)$
respectively, it suffices to consider the situation that $\Phi$ is
determined by a datum $\mathcal{A}_{\Phi}=\left(\alpha,\mu,1,0\right)$,
where $\mu\in k^{*}$ and $\mu^{s}=1$. So the relation $\sigma_{\alpha\left(i\right)}\left(x\right)=\mu\sigma_{i}\left(x\right)+c\left(x\right)$
holds for every $i=1,\ldots,r$.

2°a) If $\mu^{s}=1$ but $\mu^{s'}\neq1$ for every $s'=1,\ldots,s-1$, then,
letting $\tau\left(x\right)=\left(1-\mu\right)^{-1}c\left(x\right)$
and $\tilde{\sigma}_{i}\left(x\right)=\sigma_{i}\left(x\right)-\tau\left(x\right)$
for every $i=1,\ldots,r$, we arrive at the relation $\tilde{\sigma}_{\alpha\left(i\right)}\left(x\right)=\mu\tilde{\sigma}_{i}\left(x\right)$
for every $i=1,\ldots,r$. Furthermore, if $i_{0}$ is a unique fixed
point of $\alpha$ then $\tilde{\sigma}_{i_{0}}\left(x\right)=0$
as $\sigma_{i_{0}}\left(x\right)=\tau\left(x\right)$. So we conclude
that $P\left(x,y+\tau\left(x\right)\right)=y^{i}\tilde{P}\left(x,y^{s}\right)$
where $i=0,1$ and where $s$ denotes the length of the nontrivial
cycles in $\alpha$. The corresponding automorphism is of type (d).

2°b) If $\mu=1$ then $\alpha$ is fixed point free by virtue of Lemma
\ref{pro:auto_data} and $\textrm{char}\left(k\right)=s$, where $s$
denotes the common length's of the cycles occurring in $\alpha$. Moreover,
$s'\cdot c\left(0\right)\neq0$ for every $s'=1,\ldots,s-1$ and $\sigma_{i_{m}}\left(x\right)=\sigma_{i_{1}}\left(x\right)+\left(m-1\right)\cdot c\left(x\right)$
for every index $i_{m}$ occurring in a cycle $\left(i_{1},\ldots,i_{s}\right)$
of length $s$ in $\alpha$. Letting $r=ds$, we may suppose up to
a reordering that $\alpha$ decomposes as the product of the standard
cycles $\left(is+1,is+2\ldots,\left(i+1\right)s\right)$, where $i=0,\ldots,d-1$.
Letting $R\left(x,y\right)=\prod_{m=1}^{s}\left(y-m\cdot c\left(x\right)\right)=y^{s}-c\left(x\right)^{s-1}y$,
we conclude that \[
P\left(x,y\right)=\prod_{i=0}^{d-1}R\left(x,y-\sigma_{is}\left(x\right)\right)=\tilde{P}\left(x,y^{s}-c\left(x\right)^{s-1}y\right)\]
 for a suitable polynomial $\tilde{P}\left(x,y\right)\in k\left[x,y\right]$.
The corresponding automorphism is of type (e).

\end{proof}

\begin{enavant}
In the proof of  Theorem \ref{thm:Main_auto_thm} above, we used  the fact that every automorphism $\Phi$ of a  Danielewski surface $S=S_{\sigma,h}$, where $h\geq 2$,  is determined by a certain datum $\mathcal{A}_{\Phi}=\left(\alpha,\mu,a,b\left(x\right)\right)\in\mathfrak{S}_{r}\times k^{*}\times k^{*}\times k\left[x\right]$ for which the polynomial $\tilde{c}(x)=\sigma_{\alpha\left(i\right)}\left(ax\right)-\mu\sigma_{i}\left(x\right)\in k\left[x\right]$ does not depend on the index $i$. Actually, we needed the following more precise result.
\end{enavant}

\begin{lem}
\label{pro:auto_data} The elements in a datum $\mathcal{A}_{\Phi}=\left(\alpha,\mu,a,b\left(x\right)\right)$ corresponding to an automorphism $\Phi$ of $S$ satisfy the following additional properties

1\emph{)} The permutation $\alpha$ is either trivial or has at most
a unique fixed point. If it is nontrivial then all nontrivial cycles
with disjoint support occurring in a decomposition of $\alpha$ have
the same length $s\geq2$.

2\emph{)} If $\alpha$ is trivial then $\mu=1$ and the converse
also holds provided that $\textrm{char}\left(k\right)\neq s$. Otherwise,
if $\alpha$ is nontrivial and $\textrm{char}\left(k\right)\neq s$
then $\mu^{s}=1$ but $\mu^{s'}\neq1$ for
every $1\leq s'<s$.\\

\end{lem}
\begin{proof}
To simplify the notation, we let
$y_{i}=\sigma_{i}\left(0\right)$ for every $i=1,\ldots,r$. Note
that by hypothesis, $y_{i}\neq y_{j}$ for every $i\neq j$.
If $\alpha\in\mathfrak{S}_{r}$ has at least two fixed points, say
$i_{0}$ and $i_{1}$, then $y_{i_{0}}\left(1-\mu\right)=y_{i_{1}}\left(1-\mu\right)=\tilde{c}\left(0\right)$,
and so, $\mu=1$ and $\tilde{c}\left(0\right)=0$ as $y_{i_{0}}\neq y_{i_{1}}$. In
turn, this implies that $\alpha$ is trivial. Indeed, otherwise there
would exist an index $i$ such that $\alpha\left(i\right)\neq i$
but $y_{\alpha\left(i\right)}=y_{i}$, in contradiction with our hypothesis.
Suppose from now that $\alpha$ is nontrivial and let $s\geq2$ be the infimum of
the length's of the nontrivial cycles occurring in decomposition of
$\alpha$ into a product of cycles with disjoint supports. We deduce that $y_{i}\left(1-\mu^{s}\right)=y_{j}\left(1-\mu^{s}\right)$ for
every pair of distinct indices $i$ and $j$ in the support of a same
cycle of length $s$. Thus $\mu^{s}=1$ as $y_{i}\neq y_{j}$ for
every $i\neq j$.
If $\mu=1$ then $s'\cdot\tilde{c}\left(0\right)\neq0$ for every
$s'=1,\ldots,s-1$. Indeed, otherwise we would have $y_{\alpha^{s'}\left(i\right)}=y_{i}+s'\cdot\tilde{c}\left(0\right)=y_{i}$
for every index $i=1,\ldots,r$ which is impossible since $\alpha$
is nontrivial. In particular, $\alpha$ is fixed-point free. On the
other hand $s\cdot\tilde{c}\left(0\right)=0$ as $y_{i}=y_{\alpha^{s}\left(i\right)}=y_{i}+s\cdot\tilde{c}\left(0\right)$
for every index $i$ in the support of a cycle of length $s$ in $\alpha$.
This is possible if and only if the characteristic of the base field
$k$ is exactly $s$. We also conclude that every cycle in $\alpha$
have length $s$ for otherwise there would exist an index $i$ such
that $\alpha^{s}\left(i\right)\neq i$ but $y_{\alpha^{s}\left(i\right)}=y_{i}+s\cdot\tilde{c}\left(0\right)=y_{i}$
in contradiction with our hypothesis.

If $\mu\neq1$ then
$\mu^{s'}\neq1$ for every $s'<s$. Indeed, otherwise there would
exist an index $i$ such that $\alpha^{s'}\left(i\right)\neq i$ but
$y_{\alpha^{s'}\left(i\right)}=\mu^{s'}y_{i}+\tilde{c}\left(0\right)\sum_{p=0}^{s'-1}\mu^{p}=y_{i}$,
which is impossible. The same argument also implies that all the nontrivial
cycles in $\alpha$ have length $s$.
\end{proof}

\begin{enavant}
By combining Theorems \ref{thm:Equivalent_charac} and \ref{thm:Main_auto_thm}, we obtain the following description of the automorphisms groups  of Danielewski surfaces $S_{Q,h}$. 
\end{enavant}

\begin{thm}\label{thm:S_Q,h autos} Let $S_{Q,h}$ be the Danielewski surface in $\mathbb{A}^3_k$  defined by the equation $x^hz-Q(x,y)=0$ and let $S_{\sigma,h}$ be one of its standard forms. Then, every automorphism of $S_{Q,h}$ is of the form $\Phi^s\circ \psi \circ \Phi_s$,  where $\psi$ belongs to the subgroup $G_{\sigma,h}$ of the automorphisms group of $\mathbb{A}^3_k$ defined in Theorem \ref{thm:Main_auto_thm} and $\Phi^s$ and $\Phi_s$ are the endomorphisms of $\mathbb{A}^3_k$ defined in  \ref{rem:def phi_s}.
\end{thm}

\begin{enavant} We have seen in \ref{txt:Non_algeb_equiv} that the embeddings $i_{Q,h}$ are not rectifiable in general and so that the isomorphisms $\phi^s$ and $\phi_s$ do not extend to algebraic automorphisms of $\mathbb{A}^3_k$. Therefore, in contrast with the case of standard embeddings $i_s$ for which every automorphisms of a Danielewski surface $S\simeq S_{\sigma,h}$ arises as the restriction of an automorphism of the ambient space $\mathbb{A}^3_k$, the above result may lead one to suspect that for a general embedding $i_{Q,h}$ of $S$ as a surface $S_{Q,h}$, certain automorphisms of $S$ do not extend to algebraic automorphisms $\mathbb{A}^3_k$. In the next section we give examples of embeddings for which this phenomenon occurs. However, if $k=\mathbb{C}$, Theorem \ref{thm:Embed_analytic_equiv} leads on the contrary the following result.
\end{enavant}

\begin{cor}\label{cor:analytic extension}
Every algebraic automorphism of a Danielewski surface $S_{Q,h}$ in $\mathbb{A}^3_{\mathbb{C}}$ is extendable to a \emph{holomorphic} automorphism of $\mathbb{A}^3_{\mathbb{C}}$.
\end{cor}

\section{Special Danielewski surfaces and multiplicative group actions }

\indent\newline\noindent In this section, we fix a base field $k$ of characteristic zero and we consider special Danielewski surfaces $S$ admitting a nontrivial action of the multiplicative group $\mathbb{G}_m=\mathbb{G}_{m,k}$.  We establish that every such surface is isomorphic to a Danielewski surface $S_{Q,h}$ which admits a standard embedding in $\mathbb{A}^3_k$ as a surface defined by an equation of the form $x^hz-P\left(y\right)=0$ for a suitable polynomial $P\left(y\right)\in k\left[y\right]$. In this embedding, every multiplicative group action on $S$ arises as the restriction of a linear $\mathbb{G}_m$-action on $\mathbb{A}^3_k$. We show on the contrary that this is not the case for a general embedding of $S$ as a surface $S_{Q,h}$.

\subsection{Multiplicative group actions on special Danielewski surfaces}

\indent\newline\noindent
Every Danielewski surface isomorphic to a surface $S_{P,h}$ in $\mathbb{A}^3_k$ defined by an equation of the form $x^hz-P\left(y\right)=0$ for a certain polynomial $P\left(y\right)$ admits an nontrivial action of the multiplicative group $\mathbb{G}_m$ which arises as the restriction of the $\mathbb{G}_m$-action $\Psi$ on $\mathbb{A}_{k}^{3}$ defined by $\Psi\left(a;x,y,z\right)=H_{a}\left(x,y,z\right)=\left(ax,y,a^{-h}z\right)$. In the setting of Lemma \ref{pro:auto_data} above, the automorphisms
$H_{a}$ correspond to data $\mathcal{A}_{\phi_{a}}=\left(1,1,a,0\right)$, where $a\in k^{*}$. Here we establish that Danielewski surfaces isomorphic to a surface $S_{P,h}$ in $\mathbb{A}^3_k$ are characterised by the fact that they admit such a nontrivial $\mathbb{G}_m$-action.

\begin{enavant} \label{txt:OMaKL_surface_Tree_charac} By virtue of example \ref{exa:Main_example}
above, the collection of polynomials $\sigma_{i}\left(x\right)$,
$i=1,\ldots,r$, corresponding to a Danielewski surface $S_{P,h}\subset\mathbb{A}_{k}^{3}$
is given by $\sigma_{i}\left(x\right)=y_{i}$ for every $i=1,\ldots,r$, where $y_1,\ldots,y_r$ denote the roots of the polynomial $P$. In turn, we deduce from Theorem \ref{thm:Iso_classes} and Proposition \ref{thm:Normal_forms_iso}
above that a Danielewski surface $S_{Q,h}$ with a standard form $S_{\sigma,h}$
defined by a datum $\left(r,h,\sigma=\left\{ \sigma_{i}\left(x\right)\right\} _{i=1,\ldots,r}\right)$
is isomorphic to a surface $S_{P,h}$ as above if and only if there
exists a polynomial $\tau\left(x\right)\in k\left[x\right]$ such
that $\sigma_{i}\left(x\right)=\sigma_{i}\left(0\right)+\tau\left(x\right)$
for every $i=1,\ldots,r$. So we conclude that every such surface
correspond to a fine $k$-weighted rake $\gamma$ of the following
type.

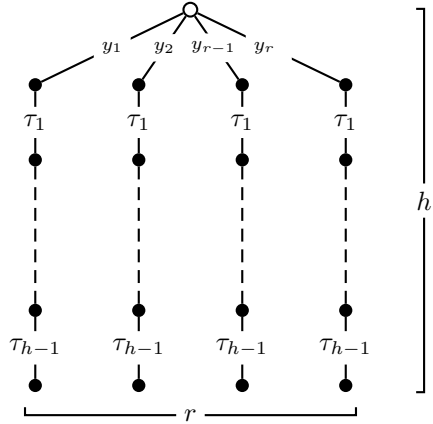
\begin{figure}[h]

\begin{pspicture}(-2.6,2.5)(8,-2.9)

\def\dedge{\ncline[linestyle=dashed]}

\rput(2,0){

\pstree[treemode=D,radius=2.5pt,treesep=1.2cm,levelsep=1cm]{\Tc{3pt}}{

  \pstree{\TC*\mput*{{\scriptsize $y_1$}}} {

              \pstree{\TC*\mput*{$\tau_1$}}{ \skiplevels{1}

                 \pstree{\TC*[edge=\dedge]}{

                    \TC*\mput*{$\tau_{h-1}$}

                                                               }

                                                      \endskiplevels

                                                     }

                                          }

\pstree{\TC*\mput*{{\scriptsize $y_2$}}}{

            \pstree{\TC*\mput*{$\tau_1$}}{\skiplevels{1}

                     \pstree{\TC*[edge=\dedge]} {

                              \TC*\mput*{$\tau_{h-1}$}

                                                                    }

                                                          \endskiplevels

                                                           }

                                                }

\pstree{\TC*\mput*{{\scriptsize $y_{r-1}\;$}}} {

              \pstree{\TC*\mput*{$\tau_1$}}{ \skiplevels{1}

                 \pstree{\TC*[edge=\dedge]}{

                    \TC*\mput*{$\tau_{h-1}$}

                                                               }

                                                      \endskiplevels

                                                     }

                                          }

\pstree{\TC*\mput*{{\scriptsize $y_r$} }}{

            \pstree{\TC*\mput*{$\tau_1$}}{\skiplevels{1}

                     \pstree{\TC*[edge=\dedge]} {

                              \TC*\mput*{$\tau_{h-1}$}

                                                                    }

                                                          \endskiplevels

                                                           }

                                                }

}

}

\pnode(-0.2,-2.8){A}\pnode(4.2,-2.8){B}

\ncbar[angleA=270, arm=3pt]{A}{B}\ncput*[npos=1.5]{$r$}

\pnode(5,2.5){C}\pnode(5,-2.6){D}

\ncbar[arm=3pt]{C}{D}\ncput*[npos=1.5]{$h$}

\end{pspicture}

\end{figure}

\end{enavant}

\begin{enavant} One can easily deduce from the description of the automorphism group  of a Danielewski surface $S_{\sigma,h}$ given Theorem \ref{thm:Main_auto_thm} above that such a surface admits a nontrivial $\mathbb{G}_m$-action if and only if it is isomorphic to a surface $S_{P,h}$. More generally, we have the following result.
\end{enavant}

\begin{thm}
\label{pro:OMaKL_surfaces_action_charac}
A special Danielewski surface $S$ admits a nontrivial action of the multiplicative group $\mathbb{G}_m$ if and only if it is isomorphic to a surface $S_{P,h}$ in $\mathbb{A}^3_k$ defined by the equation $x^hz-P(y)=0$.
\end{thm}
\begin{proof}
We may suppose that $S=S\left(\gamma\right)$ is the Danielewski surface associated with a fine $k$-weighted tree $\gamma=\left(\Gamma,w\right)$ with $r\geq 2$ elements at level $1$ and with all its leaves at level $h\geq 1$. We denote by  $\sigma=\left\{\sigma_i\left(x\right)\right\}_{i=1,\ldots,r}$ the collection of polynomial associated with $\gamma$ (see \ref{pro:WeightedTree_2_DanSurf}). By virtue of Theorem \ref{thm:Iso_classes} above, the collection $\tilde{\sigma}$ defined by \[\tilde{\sigma}_i\left(x\right)=\sigma_i\left(x\right)-\frac{1}{r}\sum_{i=1}^r\sigma_i\left(x\right)\quad i=1,\ldots,r\] leads to a Danielewski surface isomorphic to $S$. So we may suppose from the beginning that $\sigma_1\left(x\right)+\cdots + \sigma_r\left(x\right)=0$.  If $h=1$ then it follows that $S$ is isomorphic to a surface in $\mathbb{A}^3_k$ defined by an equation of the form $xz-P\left(y\right)=0$, and so, the assertion follows from the above discussion. Otherwise, if $h\geq 2$ then it follows from Theorem \ref{thm:Comb_ML_Trivial} that the structural $\mathbb{A}^1$-fibration $\pi=\pi_{\gamma}:S=S\left(\gamma\right)\rightarrow \mathbb{A}^1_k$ is unique up to automorphisms of the base. We consider $S$ as an $\mathbb{A}^1$-bundle $\rho:S\rightarrow X\left(r\right)$ defined by the transition cocycle \[g=\left\{g_{ij}=x^{-h}\left(\sigma_j\left(x\right)-\sigma_i\left(x\right)\right)\right\}_{i,j=1,\ldots,r}.\]
The same argument as in the proof of Theorem \ref{thm:Normal_forms_iso} implies that every automorphism $\Phi$ of $S$  is determined by a datum $\mathcal{A}_{\Phi}=\left(\alpha,\mu,a,b\left(x\right)\right)\in\mathfrak{S}_{r}\times k^{*}\times k^{*}\times k\left[x\right]$ for which the polynomial $\sigma_{\alpha\left(i\right)}\left(ax\right)-\mu\sigma_{i}\left(x\right)\in k\left[x\right]$ does not depend on the index $i$. In view of the composition rule given in the same proof, we deduce that an automorphism $\Phi$ of $S$ may belong to a subgroup of $\textrm{Aut}\left(S\right)$ isomorphic to $\mathbb{G}_m$ only if its associated datum is of the form $\mathcal{A}_{\Phi}=\left(\alpha,\mu,a,0\right)$. Suppose that there exists a nontrivial automorphism  $\Phi$ determined by such a datum $\mathcal{A}_{\Phi}$. Then, since $\alpha\in\mathfrak{S}_r$, there exists an integer $N\geq 1$ such that
the polynomial $c\left(x\right)=\sigma_i\left(a^Nx\right)-\mu^N\sigma_i\left(x\right)$ does not depend on the index $i=1,\ldots,r$. Since $\sigma_1\left(x\right)+\cdots+\sigma_r\left(x\right)=0$ by hypothesis, we conclude that the identity $\sigma_i\left(a^Nx\right)=\mu^N\sigma_i\left(x\right)$ holds for every index $i=1,\ldots,r$. In particular, it follows that $\sigma_i\left(0\right)=\mu^N\sigma_i\left(0\right)$ for every index $i=1,\ldots,r$. Thus $\mu^N=1$ since $\gamma$ is a fine $k$-weighted tree with at least two elements at level $1$. Suppose that one of the polynomials $\sigma_i$ is not constant. Then the above identity implies that $a^{Np}=1$ for a certain integer $p$. Therefore, every automorphism $\Phi$ of $S$ with associated datum $\left(\alpha,\mu,a,0\right)$ is cyclic and $\textrm{Aut}\left(S\right)$ can not contain a subgroup isomorphic to $\mathbb{G}_m$. So, $S$ admits a nontrivial $\mathbb{G}_m$-action only if the polynomials $\sigma_i$, $i=1,\ldots,r$ are constant.  This completes the proof since these fine $k$-weighted trees correspond to Danielewski surfaces $S_{P,h}$ by virtue of \ref{txt:OMaKL_surface_Tree_charac} above.
\end{proof}

\subsection{ Extensions of multiplicative group actions on a Danielewski surface}

 \indent\smallskip\newline\noindent It follows from Theorem \ref{thm:Main_auto_thm}  that every special Danielewski surface $S$ equipped with a nontrivial $\mathbb{G}_m$-action admits an equivariant embedding in $\mathbb{A}^3_k$ as a surface $S_{P,h}$ defined by an equation of the form $x^hz-P\left(y\right)=0$. In this embedding, the $\mathbb{G}_m$-action on $S$ even arises as the restriction of a linear $\mathbb{G}_m$-action on $\mathbb{A}^3_k$ corresponding to automorphisms of type b) in \ref{thm:Main_auto_thm}. On the other hand, a surface $S$ isomorphic to a surface $S_{P,h}$ admits closed embeddings $i_{Q,h}:S\hookrightarrow \mathbb{A}^3_k$ in $\mathbb{A}^3_k$ as surfaces $S_{Q,h}$ defined by equations of the form $x^hz-R\left(x,y\right)P\left(y\right)=0$ (see Theorem \ref{thm:Equivalent_charac}). It is natural to ask if there always exists $\mathbb{G}_m$-actions on $\mathbb{A}^3_k$ making these general embeddings equivariant. Clearly, this holds if the embedding $i_{Q,h}$ is algebraically equivalent to a standard embedding of $S$ as a surface $S_{P,h}$. The following result shows that there exists non rectifiable closed embeddings $i_{Q,h}$ of $S$  for which no nontrivial $\mathbb{G}_m$-action on $S$ can be extended to an action on the ambient space.

\begin{thm}\label{txt:Multiplicative_action} 
\label{thm:Multiplicative_auto_non_extension} Every Danielewski surface $S\subset\mathbb{A}_{k}^{3}$
defined by the equation $x^{h}z-\left(1-x\right)P\left(y\right)=0$, where $h\geq2$
and where $P\left(y\right)$ has $r\geq2$ simple roots, admits a nontrivial $\mathbb{G}_m$-action
$\tilde{\theta}:\mathbb{G}_m\times S\rightarrow S$ which is not
algebraically extendable to $\mathbb{A}_{k}^{3}$. More precisely,
for every $a\in k\setminus\left\{ 0,1\right\} $ the automorphism
$\tilde{\theta}_{a}=\tilde{\theta}\left(a,\cdot\right)$ of $S$ do
not extend to an algebraic automorphism of $\mathbb{A}_{k}^{3}$.
\end{thm}

\begin{proof}
The endomorphisms  $\Phi^s$ and $\Phi_s$  of $\mathbb{A}^3_k$ defined by $\Phi^s\left(x,y,z\right)=\left(x,y,(1-x\right)z)$
and $\Phi_s\left(x,y,z\right)=\left(x,y,(\sum_{i=0}^{h-1}x^{i})z+P\left(y\right)\right)$
induce isomorphisms $\phi^s$ and $\phi_s$ between $S$ and the
surface $S_{P,h}$ defined by the equation $x^{h}z-P\left(y\right)=0$ (see  \ref{rem:def phi_s}). The
latter admits an action $\theta:\mathbb{G}_m\times S_{P,h}\rightarrow S_{P,h}$
of the multiplicative group $\mathbb{G}_m$ defined by $\theta\left(a,x,y,z\right)=H_{a}\left(x,y,z\right)=\left(ax,y,a^{-h}z\right)$  for every $a\in k^{*}$. The corresponding action $\tilde{\theta}$
on $S$ is therefore defined by $\tilde{\theta}\left(a,x,y,z\right)=\tilde{\theta}_{a}\left(x,y,z\right)=\phi^s\circ H_{a}\left(x,y,z\right)\mid_{S_{P,h}}\circ\phi_s$. 
Since by construction, $\tilde{\theta}_{a}^{*}\left(x\right)=ax$
for every $a\in k^{*}$, the assertion is a consequence of the following
Lemma which guarantees that the automorphisms $\tilde{\theta}_{a}$
of $S$ are not algebraically extendable to an automorphism of $\mathbb{A}_{k}^{3}$ for every $a\in k^*\setminus \{ 1 \}$.
\end{proof}

\begin{lem}
\label{lem:Rigid_extension} If $\Phi$ is an algebraic automorphism
of $\mathbb{A}_{k}^{3}$ extending an automorphism of $S$, then $\Phi^{*}\left(x\right)=x$.
\end{lem}
\begin{proof}
Our proof is similar to the one of Theorem 2.1 in \cite{MoP05}. We
let $\Phi$ be an automorphism of $\mathbb{A}_{k}^{3}$ extending
an arbitrary automorphism of $S$. Since $f_{1}=x^{h}z-\left(1-x\right)P\left(y\right)$
is an irreducible polynomial, there exists $\mu\in k^{*}$ such that
$\Phi^{*}\left(f_{1}\right)=\mu f_{1}$. Therefore, for every $t\in k$,
the automorphism $\Phi$ induces an isomorphism between the level
surfaces $f_{1}^{-1}\left(t\right)$ and $f_{1}^{-1}\left(\mu^{-1}t\right)$
of $f_{1}$. There exists an open subset $U\subset\mathbb{A}_{k}^{1}$
such that for every $t\in U$, $f_{1}^{-1}\left(t\right)$ is a special
Danielewski surfaces isomorphic to a one defined by a fine $k$-weighted
rake $\gamma$ whose underlying tree $\Gamma$ is isomorphic to the
one associated with $S$. Since $\Gamma$ is not a comb, it follows
from Theorem \ref{thm:Comb_ML_Trivial} that for every $t\in U$,
the projection ${\rm pr}_{x}:f_{1}^{-1}\left(t\right)\rightarrow\mathbb{A}_{\mathbb{C}}^{1}$
is a unique $\mathbb{A}^{1}$-fibration on $f_{1}^{-1}\left(t\right)$
up to automorphisms of the base. Furthermore, ${\rm pr}_{x}:f_{1}^{-1}\left(t\right)\rightarrow\mathbb{A}_{k}^{1}$
has a unique degenerate fiber, namely ${\rm pr}_{x}^{-1}\left(0\right)$.
Therefore, for every $t\in U$, the image of the ideal $\left(x,f_{1}-t\right)$
of $k\left[x,y,z\right]$ by $\Phi^{*}$ is contained in the ideal
$\left(x,\mu f_{1}-t\right)=\left(x,P\left(y\right)+\mu^{-1}t\right)$,
and so $\Phi^{*}\left(x\right)\in\bigcap_{t\in U}\left(x,P\left(y\right)+\mu^{-1}t\right)=\left(x\right)$.
Since $\Phi$ is an automorphism of $\mathbb{A}_{k}^{3}$, we conclude
that there exists $c\in k^{*}$ such that $\Phi^{*}\left(x\right)=cx$.
In turn, this implies that for every $t,u\in k$, $\Phi$ induces
an isomorphism between the surfaces $S_{t,u}$ and $\tilde{S}_{t,u}$
defined by the equations $f_{1}+tx+u=x^{h}z-\left(1-x\right)P\left(y\right)+tx+u=0$
and $f_{1}+\mu^{-1}ctx+\mu^{-1}u=x^{h}z-\left(1-x\right)P\left(y\right)+\mu^{-1}ctx+\mu^{-1}u=0$
respectively. Since $\deg\left(P\right)\geq2$ there exists $y_{0}\in k$
such that $P'\left(y_{0}\right)=0$. Note that $y_{0}$ is not a root of $P$ as these ones are simple. We let $t=-u=-P\left(y_{0}\right)$.
 Since $h\geq2$, it follows from the Jacobian Criterion that $S_{t,u}$
is singular, and even non normal along the nonreduced component of
the fiber ${\rm pr}_{x}^{-1}\left(0\right)$ defined by the equation $\left\{ x=0;\, y=y_{0}\right\} $.
Therefore $\tilde{S}_{t,u}$ must be singular along a multiple component
of the fiber ${\rm pr}_{x}^{-1}\left(0\right)$. This the case if
and only if the polynomial $P\left(y\right)-\mu^{-1}cP\left(y_{0}\right)$
has a multiple root, say $y_{1}$, such that $P\left(y_{1}\right)-\mu^{-1}P\left(y_{0}\right)=0$.
Since $P\left(y_{0}\right)\neq0$ this condition is satisfied if and
only if $c=1$. This completes the proof.
\end{proof}
\begin{example}
In particular, even the involution of the surface $S$ defined by the equation $x^{2}z-\left(1-x\right)P\left(y\right)=0$ induced by the endomorphism
$J\left(x,y,z\right)=\left(-x,y,\left(1+x\right)\left(\left(1+x\right)z+P\left(y\right)\right)\right)$
of $\mathbb{A}^3_k$ does not extend to an algebraic automorphism of $\mathbb{A}_{k}^{3}$.
\end{example}

\noindent It turns out that this kind of phenomenon does not occur with additive group actions. More precisely, we have the following result.

\begin{prop}
Let $S_{Q,h}$ be the Danielewski surface in $\mathbb{A}_{k}^{3}$ defined by the equation $x^hz-Q(x,y)=0$.
Then, every $\mathbb{G}_a$-action on $S_{Q,h}$ arrises as the restriction of a $\mathbb{G}_a$-action on $\mathbb{A}_{k}^{3}$ defined by 
$\tilde{\Delta}\left(t,x,y,z\right)=\left(x, y+x^hb(x)t, z+x^{-h}(Q(x,y+x^hb(x)t)-Q(x,y))\right))$,
for a certain polynomial $b\left(x\right)\in k[x]$.
\end{prop}

\begin{proof}
With the notation of  Remark \ref{rem: k+-actions}, it follows from  Theorem \ref{thm:S_Q,h autos} that every additive group action on $S_{Q,h}$ is induced by the restriction to $S_{Q,h}$ of a collection of endomorphisms of $\mathbb{A}^3_k$ of the form $\delta_{t,b}=\Phi^s\circ \Delta_{tb(x)} \circ \Phi_s$, where $b\in k\left[x\right]$. One checks that  
\[\begin{array}{lcl} \delta_{t,b}(x,y,z) &=& (x, y+x^hb(x)t, z+x^{-h}(Q(x,y+x^hb(x)t)-Q(x,y))+\alpha(x,y)(x^hz-Q(x,y))), \end{array}\] 
\noindent for a certain polynomial $\alpha(x,y)\in k[x,y]$. Note that if $\alpha\left(x,y\right)\neq 0$, these endomorphisms $\delta_{t,b}$ do not define a $\mathbb{G}_a$-action on $\mathbb{A}^3_k$. However, they induce an action on $S_{Q,h}$ which coincides with the one induced by the $\mathbb{G}_a$-action $\tilde{\Delta}$ above. 
\end{proof}

\begin{enavant} If $k=\mathbb{C}$, Corollary \ref{cor:analytic extension} implies in particular that every automorphism of $S$  extends to an holomorphic
automorphism of $\mathbb{A}_{\mathbb{C}}^{3}$. This leads the following result which contrasts
with an example, given by H. Derksen, F. Kutzschebauch and J. Winkelmann
in \cite{DKW}, of a non-extendable $\mathbb{C}_{+}$-action on an
hypersurface in $\mathbb{A}_{\mathbb{C}}^{5}$ which is even holomorphically
inextendable .
\end{enavant}
\begin{prop}
Every surface $S\subset\mathbb{A}_{\mathbb{C}}^{3}$ defined by the equation
$x^{h}z-\left(1-x\right)P\left(y\right)=0$, where $h\geq2$ and where
$P\left(y\right)$ has $r\geq2$ simple roots, admits a nontrivial
$\mathbb{C}^{*}$-action which is algebraically inextendable but holomorphically
extendable to $\mathbb{A}_{\mathbb{C}}^{3}$.
\end{prop}
\begin{proof}
We let $\tilde{\theta}:\mathbb{C}^{*}\times S\rightarrow S$ be the
$\mathbb{C}^{*}$-action on the surface $S\subset\mathbb{A}_{\mathbb{C}}^{3}$
defined by the equation $x^{2}z-\left(1-x\right)P\left(y\right)=0$ constructed
in the proof of Theorem \ref{txt:Multiplicative_action}.
For every $a\in\mathbb{C}^{*}$, the automorphism $\tilde{\theta}\left(a,\cdot\right)$
of $S$ maps a closed point $\left(x,y,z\right)\in S$ to the point
$\tilde{\theta}\left(a,x,y,z\right)=\left(ax,y,a^{-2}\left(1-ax\right)\left(\left(1+x\right)z+P\left(y\right)\right)\right)$.
One checks that the holomorphic automorphism $\Phi_{a}$ of $\mathbb{A}_{\mathbb{C}}^{3}$
such that $\Phi_{a}\mid_{S}=\tilde{\theta}\left(a,\cdot\right)$ is
the following one: \begin{eqnarray*}
\Phi_{a}\left(x,y,z\right) & = & \left(ax,y,a^{-2}e^{\left(1-a\right)x}z+\left(ax\right)^{-2}P\left(y\right)\left(e^{\left(1-a\right)x}\left(x-1\right)-ax+1\right)\right).\end{eqnarray*}
 Clearly, the holomorphic map $\Phi:\mathbb{C}^{*}\times\mathbb{A}_{\mathbb{C}}^{3}\rightarrow\mathbb{A}_{\mathbb{C}}^{3}$,
$\left(a,\left(x,y,z\right)\right)\mapsto\Phi_{a}\left(x,y,z\right)$
defines a $\mathbb{C}^{*}$-action on $\mathbb{A}_{\mathbb{C}}^{3}$
extending the one $\tilde{\theta}$ on $S$.
\end{proof}

\bibliographystyle{amsplain}

\end{document}